\documentclass[11pt,reqno]{amsart}
\usepackage{amsmath,amssymb,amsthm,mathtools}
\usepackage[margin=1.1in]{geometry}
\usepackage{enumitem}
\usepackage{graphicx,booktabs,tikz}
\usepackage{algorithm}
\usepackage{algpseudocode}
\graphicspath{{figs/}}
\usepackage{xcolor}
\usepackage[colorlinks=true,linkcolor=blue!55!black,citecolor=blue!55!black,urlcolor=blue!55!black]{hyperref}
\usepackage{cite}

\numberwithin{equation}{section}
\newtheorem{theorem}{Theorem}[section]
\newtheorem{proposition}[theorem]{Proposition}
\newtheorem{lemma}[theorem]{Lemma}
\newtheorem{corollary}[theorem]{Corollary}
\newtheorem{problem}[theorem]{Problem}
\theoremstyle{definition}
\newtheorem{definition}[theorem]{Definition}
\newtheorem{assumption}[theorem]{Assumption}
\theoremstyle{remark}
\newtheorem{remark}[theorem]{Remark}

\newcommand{\R}{\mathbb{R}}
\newcommand{\dr}{\,dr}
\newcommand{\dx}{\,dx}
\newcommand{\dS}{\,dS}
\newcommand{\dt}{\,dt}

\newcommand{\norm}[1]{\lVert #1\rVert}
\newcommand{\abs}[1]{\lvert #1\rvert}
\newcommand{\pair}[1]{\langle #1\rangle}
\newcommand{\Om}{\Omega}
\newcommand{\bv}{\mathbf{v}}
\newcommand{\bj}{\mathbf{j}}   
\newcommand{\bnu}{\boldsymbol{\nu}}
\newcommand{\GD}{\partial\Omega_D}
\newcommand{\GN}{\partial\Omega_N}
\newcommand{\HD}{H^1_D(\Om)}
\newcommand{\HDd}{H^{-1}_D(\Om)}
\newcommand{\ee}{\mathcal{E}}
\newcommand{\tp}{\tilde p}
\newcommand{\tn}{\tilde n}

\newcommand{\Y}{\mathcal{Y}}
\newcommand{\Rdata}{\mathcal{R}_{\mathrm{data}}}
\newcommand{\RF}{\mathcal{R}_{F}}
\newcommand{\Jou}{\mathcal{J}}
\DeclareMathOperator{\esssup}{ess\,sup}

\begin{document}

\title[Global existence for a MEMS discharge model]{Global existence for a drift--diffusion model of electrical discharge in MEMS}


\author{Runan He}
\address{Instituto de Ciencias Matem\'aticas (ICMAT), Calle Nicol\'as Cabrera 13, 28049 Madrid, Spain}

\subjclass[2020]{35K57, 35Q81, 35B45, 35M33, 82D10}
\keywords{Drift--diffusion system, gas discharge, MEMS, Townsend ionization, global existence, blow-up, continuation criterion, electrostatic energy}

\begin{abstract}
We study a coupled elliptic--parabolic system modelling electrical discharge in a micro-electro-mechanical (MEMS) capacitor: a Poisson equation for the electrostatic potential, coupled to drift--diffusion equations for the ion and electron densities with a Townsend ionization source whose rate per electron grows linearly in the electric field. Gimperlein, He and Lacey proved local well-posedness and left open whether the densities blow up in finite time. We first prove a continuation criterion in a weak norm: a finite-time singularity forces the space charge to be unbounded in $H^{-1}_D$, the dual of the space of $H^1$ functions vanishing on the electrodes, so blow-up cannot be quasi-neutral. We then exclude it altogether when the two diffusivities coincide and the vacuum field is bounded: the weak solution is then global in time. To our knowledge this is the first global existence result without smallness assumptions for a
drift--diffusion--Poisson model of gas discharge whose ionization rate is unbounded in the field. The proof uses the electrostatic energy and the number of carriers in place of the usual free energy: the ionization source cancels in the energy balance, the Joule heating is dissipative, and the two quantities close through a linear Gronwall inequality. Numerical experiments on a discharge in a micrometre gap, with the Townsend constants and the mobilities of air, provide evidence for the extension to unequal diffusivities and bear on the existence of a breakdown threshold.
\end{abstract}

\maketitle

\section{Introduction}

\subsection{The model}\label{ssec:model}
We are concerned with the following model for electrical discharge in the gas-filled gap of a MEMS capacitor, which was proposed by Gimperlein, He and Lacey in \cite{GHL26}, in the spirit of Budd's corona model \cite{Bud91}:
\begin{subequations}\label{eq:sys}
\begin{align}
-\Delta\phi &= \frac{p-n}{\epsilon_0} &&\text{in }\Om\times(0,T),\label{eq:poisson}\\
\partial_t p-\nabla\cdot\big(\epsilon_+\nabla p+p\,(\mu_+\nabla\phi-\bv)\big)&=F(n,|\nabla\phi|) &&\text{in }\Om\times(0,T),\label{eq:p}\\
\partial_t n-\nabla\cdot\big(\epsilon_-\nabla n-n\,(\mu_-\nabla\phi+\bv)\big)&=F(n,|\nabla\phi|) &&\text{in }\Om\times(0,T),\label{eq:n}
\end{align}
\end{subequations}
with the ionization source
\begin{equation}\label{eq:F}
F(n,|\nabla\phi|)=\Gamma(|\nabla\phi|)\,n,\qquad
\Gamma(E):=\mu_-E\big(\alpha_1e^{-\alpha_2/E}-\eta_0\big)\ \ (E>0),\qquad \Gamma(0):=0,
\end{equation}
subject to the initial conditions
\begin{equation}\label{eq:ic}
p(\cdot,0)=p_0,\qquad n(\cdot,0)=n_0\qquad\text{in }\Om,
\end{equation}
and the mixed boundary conditions
\begin{equation}\label{eq:bc}
\begin{aligned}
\phi&=\phi_D, & p&=p_D, & n&=n_D && \text{on }\GD\times(0,T),\\
\bnu\cdot\nabla\phi&=0, & \bnu\cdot\nabla p&=0, & \bnu\cdot\nabla n&=0 && \text{on }\GN\times(0,T).
\end{aligned}
\end{equation}
Here $\phi$ is the electrostatic potential, and $p$ and $n$ are the densities of positive ions and of electrons. The bounded domain $\Om\subset\R^d$, $d\in\{2,3\}$, represents the gas-filled part of the interior of a MEMS device; its boundary is the disjoint union of a Dirichlet part $\GD$ and a Neumann part $\GN$, and $\bnu$ is the outward unit normal. The Dirichlet part $\GD$ consists of the two electrodes, on which the potential and the carrier densities are prescribed. The Neumann part $\GN$ is an artificial boundary that closes the computational region away from the gap. 

The constant $\epsilon_0>0$ is the permittivity, $\epsilon_\pm>0$ are diffusivities, $\mu_\pm>0$ are mobilities, and $\alpha_1,\alpha_2,\eta_0>0$ are the constants of the Townsend ionization law. The argument $E$ of $\Gamma$ is the field strength $|\nabla\phi|$, the factor $\alpha_1e^{-\alpha_2/E}-\eta_0$ is the net Townsend coefficient, the number of ionizing collisions less attachments to neutral molecules per electron per unit \emph{length} of drift, and $\mu_-E$ is the electron drift speed, so that $\Gamma(E)$ counts ionization events per electron per unit \emph{time}. The velocity field $\bv$ of the gas is prescribed, incompressible and tangential to the boundary: $\nabla\cdot\bv=0$ in $\Om$, $\bnu\cdot\bv=0$ on $\partial\Om$. The terms with $\epsilon_\pm$, $\mu_\pm$ and $\bv$ in \eqref{eq:p} and \eqref{eq:n} are the diffusion, drift and convection terms.

In contrast to the classical drift--diffusion systems, for example the van Roosbroeck system \cite{MRS90}, the generation term \eqref{eq:F} of \eqref{eq:sys} is linear in the electron density but has a rate that is unbounded in the field:
\begin{equation}\label{eq:Gamma-bound}
\abs{\Gamma(E)}\le\mu_-\gamma_*E\quad\text{for all }E\ge0,\qquad \gamma_*:=\max\{\alpha_1,\ \eta_0\},
\end{equation}
and, for $\alpha_1>\eta_0$, $\Gamma(E)\sim\mu_-(\alpha_1-\eta_0)E$ as $E\to\infty$, corresponding to the Townsend avalanche: the electron population grows at a rate proportional to the strength of the field that accelerates it.

The physical background and the related literature are discussed in Section~\ref{ssec:background}; the derivation of the model is given in \cite[\S1.1]{GHL26}, and the stationary problem is studied in \cite{HX26}. Under natural regularity and compatibility assumptions on the domain and the data, stated, in the form used here, in Assumption~\ref{ass:main} below, Gimperlein, He and Lacey proved that problem \eqref{eq:sys}--\eqref{eq:bc} has a unique weak solution, with $p,n\ge0$, on a time interval $(0,T_1)$ \cite[Theorem~1.1]{GHL26}; whether $T_1$ may be taken arbitrarily large, or whether instead $\norm{p(t)}_{L^2}$ and $\norm{n(t)}_{L^2}$ blow up as $t\rightarrow T_1$ for finite $T_1$, is left open there \cite[Remark~1.2]{GHL26}.

The finite existence time results from the $L^2$ energy estimate for the carrier densities, which is coupled to the field through the Poisson equation. For $Y(t):=\norm{p(t)}_{L^2}^2+\norm{n(t)}_{L^2}^2$, the $L^2$ estimate \cite[(4.14)]{GHL26}, proved there for the solutions $(p^M,n^M)$ of the truncated problem of \cite[\S3]{GHL26}, which approximate the weak solution, and written here without the superscript $M$, is an integral inequality of Riccati type:
\begin{equation}\label{eq:GHL-riccati}
Y(t)+\int_0^t\Big(\norm{\nabla p}_{L^2}^2+\norm{\nabla n}_{L^2}^2\Big)\dr
\le \tilde H_4+\tilde H_5\,t+\tilde H_6\int_0^t\big(Y+Y^2\big)\dr,
\end{equation}
which, by the nonlinear Gronwall inequality of Bihari \cite{Bih56}, bounds $Y$ only up to the blow-up time of the comparison equation (see \cite[(4.15)]{GHL26}); $T_1$ is defined accordingly, and beyond $T_1$ the estimate gives no information, so that both global existence and finite-time blow-up remain possible. The quadratic term comes from the ionization term, although the source \eqref{eq:F} is linear in $n$; Remark~\ref{rem:origin} explains how it arises.

\subsection{Main results}
The purpose of this paper is to determine whether the finite existence time is a feature of the estimate \eqref{eq:GHL-riccati} or of the system. We prove a continuation criterion for arbitrary diffusivities, and global existence when the two diffusivities coincide and the vacuum field is bounded. Our results are stated for the weak solutions of \cite{GHL26} (Definition~\ref{def:weak}), extended to their maximal interval of existence $(0,T_{\max})$ (Section~\ref{ssec:local}). Throughout, $\HDd$ denotes the dual of $\HD:=\{u\in H^1(\Om):\ u=0\ \text{on }\GD\}$, normed by $\norm{\nabla u}_{L^2}$ (Section~\ref{ssec:notation}), and
\[
\kappa:=1\ \text{ if }d=2,\qquad \kappa:=2\ \text{ if }d=3 .
\]

\begin{theorem}[Continuation criterion]\label{thm:A}
Let Assumption~\ref{ass:main}\,\textup{(A1)--(A4)} hold and let $(\phi,p,n)$ be the maximal weak solution of \eqref{eq:sys}--\eqref{eq:bc} on $(0,T_{\max})$. If $T_{\max}<\infty$, then
\begin{equation}\label{thm:A_1st_conclusion}
\int_0^{T_{\max}}\norm{\nabla\phi(t)}_{L^2(\Om)}^{2\kappa}\dt=\infty,\end{equation}
and consequently
\begin{equation}\label{thm:A_2nd_conclusion}\int_0^{T_{\max}}\norm{p(t)-n(t)}_{\HDd}^{2\kappa}\dt=\infty .
\end{equation}
In particular, \begin{equation}\label{thm:A_3rd_conclusion}\limsup_{t\rightarrow T_{\max}}\norm{p(t)-n(t)}_{\HDd}=\infty.\end{equation}
\end{theorem}

\begin{theorem}[Global existence]\label{thm:B}
Let Assumption~\ref{ass:main}\,\textup{(A1)--(A4)} hold, and assume in addition
\begin{itemize}[leftmargin=2.2em]
\item[\textup{(H1)}] equal diffusivities: $\epsilon_+=\epsilon_-=:\epsilon$;
\item[\textup{(H2)}] bounded vacuum field: the harmonic extension $\phi_D$ of the electrode potential in (A2) satisfies $\nabla\phi_D\in L^\infty(\Om)$.
\end{itemize}
Then $T_{\max}=\infty$: the weak solution of \cite[Theorem~1.1]{GHL26} is global in time. Moreover, there is a positive constant
\begin{multline*}
C=C\Big(d,~\Om,~ \epsilon,~\mu_\pm,~\epsilon_0,~\max\{\alpha_1,~\eta_0\},~\norm{\phi_D}_{H^2},~\norm{\nabla\phi_D}_{L^\infty},\\
~\norm{p_D}_{W^{2,\sigma}},~\norm{n_D}_{W^{2,\sigma}},~\norm{\bv}_{L^\infty((0,\infty),H^2)}\Big)
\end{multline*}
with $\sigma$ as in (A3), such that
\begin{equation}\label{eq:growth}
\begin{gathered}
\frac{\epsilon_0}{2}\norm{\nabla(\phi-\phi_D)(t)}_{L^2}^2+\int_\Om\big(p(t)+n(t)\big)\dx\le C\,N_0\,(1+t)\,e^{Ct},\\
\norm{p(t)}_{L^2}^2+\norm{n(t)}_{L^2}^2+\int_0^t\big(\norm{\nabla p}_{L^2}^2+\norm{\nabla n}_{L^2}^2\big)\dr\le C\,N_0\exp\left(C\,N_0^{\kappa}\,e^{Ct}\right)
\end{gathered}
\end{equation}
hold for all $t\ge0$,
where $N_0:=1+\norm{p_0}_{L^2}^2+\norm{n_0}_{L^2}^2$.
\end{theorem}

Theorem~\ref{thm:A} requires neither (H1) nor (H2): it holds for arbitrary diffusivities and for every electrode potential admitted by (A2). By \eqref{thm:A_3rd_conclusion}, a finite-time singularity, if it occurs, is carried by a macroscopic separation of charge. Theorem~\ref{thm:B} answers the question raised in \cite[Remark~1.2]{GHL26} for equal diffusivities and a bounded vacuum field: no finite-time blow-up occurs. Hypothesis (H1) restricts the diffusivities only; the mobilities are arbitrary. It does not hold for air, in which the electrons diffuse four orders of magnitude faster than the ions (Table~\ref{tab:params}). Whether Theorem~\ref{thm:B} extends to unequal diffusivities is open; Proposition~\ref{prop:conditional} reduces that case, under (H2), to a single a priori bound, which Section~\ref{ssec:num-uneq} evaluates on a computed solution and Section~\ref{ssec:unequal} discusses. Hypothesis (H2) is not implied by (A1)--(A2); it holds when $\GN$ consists of planar pieces that do not meet one another and the electrodes are $C^{1,1}$ surfaces that carry constant potentials and meet $\GN$ at right angles, which covers the geometry of \cite{GHL26} as long as the electrodes are apart, and it fails at touch-down (Section~\ref{ssec:touchdown}). The bounds \eqref{eq:growth} are not expected to be sharp: they allow the $L^2$-norms of the densities to grow double-exponentially in time, whereas in the simulations of Section~\ref{sec:numerics} the number of carriers grows by orders of magnitude and then levels off or decreases (Figure~\ref{fig:voltages}(a)). Section~\ref{sec:remarks} collects the open problems raised by the proofs of Theorems~\ref{thm:A} and~\ref{thm:B} and by the simulations: global existence for unequal diffusivities, the failure of (H2) at touch-down, the singularities that the two theorems do not exclude, and the breakdown threshold. Since Theorem~\ref{thm:B} excludes finite-time blow-up under (H1)--(H2), such a threshold cannot separate global existence from blow-up; we expect it to separate relaxation to a stationary state from unbounded growth as $t\to\infty$ (Problem~\ref{conj:threshold}).

\subsection{Background and related work}\label{ssec:background}
\emph{Electrostatic MEMS and touch-down.}
In an electrostatically actuated MEMS device, a voltage $V$ pulls a deformable electrode towards a fixed one against elastic restoring forces; above the pull-in voltage no equilibrium exists, and the electrodes touch in finite time (touch-down). In the vanishing-aspect-ratio limit the deflection $u$ satisfies a semilinear evolution equation with the singular nonlinearity $\lambda/(1+u)^2$ of the electrostatic force, $\lambda\propto V^2$, and touch-down is a quenching singularity; see the monographs \cite{EGG10,PB02} and the survey \cite{LW17}, and, for the stationary problem, the touch-down behaviour, quenching and hyperbolic versions, \cite{DZ19,FG93,GG07,GHL22,GH00,Guo08,Guo10,KLNT11,KLNT15,Kaw96}. Free-boundary models, in which the potential is harmonic in the deformed gap, were studied in \cite{ELW13,ELW14,LW13}, post-touch-down configurations of a regularized model in \cite{LLG14}, and models coupling the deflection to the gas pressure in the gap in \cite{GHL23,GHL24a,GHL24b}. All these models neglect space charge: the field is the vacuum field. As the gap closes, this field, of order $V/w$ with $w$ the gap width, grows large enough for electrons to multiply by impact ionization; a discharge can then form, modify the field by its space charge, and develop into a spark. We study the discharge before contact, decoupled from the mechanics: as in \cite{GHL26}, the position of the electrodes is prescribed.

\emph{Electrical breakdown in micrometre gaps.}
For a fixed gap between parallel electrodes, Townsend's criterion states that a discharge becomes self-sustained when the avalanche started by one electron produces ions that release on average at least one new electron from the cathode by secondary emission. The criterion neglects the space charge: the carrier densities obey linear equations in the vacuum field, and the breakdown voltage depends on pressure and gap width only through their product (Paschen's law) \cite{Rai91}. Below roughly ten micrometres in air at atmospheric pressure, the regime of MEMS, the measured breakdown voltage falls below the Paschen curve \cite{CYW06,DTD00,Ger59,OSE00,ST02,TD99}, which is attributed to field emission from the cathode enhanced by the approaching ions \cite{BK55,GP10,TG10}; for reviews see \cite{FZVW20,GLDV20,GV14,SB16}. Problem \eqref{eq:sys}--\eqref{eq:bc} couples the space charge to the field but includes neither secondary nor field emission, since the carrier densities on the electrodes are prescribed; the breakdown threshold of Problem~\ref{conj:threshold}, the counterpart of the breakdown voltage in this model, is therefore not given by Townsend's criterion.

\emph{Continuum models of discharges and their analysis.}
Drift--diffusion--Poisson systems with a Townsend ionization source, such as \eqref{eq:sys}, are the standard fluid description of coronas and streamers \cite{Rai91}. System \eqref{eq:sys} goes back to Budd and Wheeler \cite{Bud91,BW90}, who formulated it for a positive corona between coaxial cylinders and observed numerically that the stationary solution loses stability at high voltage. Streamers have been simulated with such models \cite{ML97} and analysed as a pattern-formation problem \cite{Ebe06}: planar ionization fronts are pulled fronts, whose width is set by diffusion and tends to zero with the diffusivity \cite{EvSC96,EvSC97}, they are unstable to transverse perturbations \cite{DEM08}, and streamers branch \cite{ArrEH02,MRE04}.

Rigorous results concern three neighbouring systems: the van Roosbroeck system of semiconductor theory, a gas-discharge model with a bounded ionization rate, and the model of Degond, Lucquin-Desreux and Morrow for a Townsend discharge. For $F=0$ and $\bv=0$, \eqref{eq:sys} is the van Roosbroeck system without doping and recombination \cite{Jue09,Mar86,MRS90}, for which existence of stationary solutions, global well-posedness and convergence to equilibrium go back to \cite{Gaj85,GG86,Moc72,Moc74,Moc75}. The Lyapunov functional in \cite{Gaj85,GG86}, as in \cite{BD00} for the Nernst--Planck--Poisson system, is the free energy of Remark~\ref{rem:freeenergy}; it is bounded from below because the electrostatic energy enters with a positive sign: like charges repel. Global existence, uniqueness and regularity were proved in \cite{Isl03,Isl06} for a gas-discharge model with the source $(\alpha-\beta p)\,n+\eta$ in both carrier equations, in which the ionization frequency $\alpha$, the recombination coefficient $\beta$ and the external ionization rate $\eta$ may depend on position, time and the field, but $\alpha$ is bounded in the field, and with nonlinear electrode conditions that couple the two densities.

The generation rate per carrier is bounded in the field in these global results, but unbounded in it for the avalanche generation $\alpha_+|J_+|+\alpha_-|J_-|$ of semiconductor theory, where the ionization coefficients $\alpha_\pm$ are bounded functions of the field and $J_\pm$ are the current densities. The source \eqref{eq:F} is of this type: it is the electron term $\alpha_-|J_-|$ with $J_-$ replaced by the drift flux $\mu_-n\nabla\phi$ and $\alpha_-$ by the net Townsend coefficient, so that $F$ grows quadratically in $(n,\nabla\phi)$ jointly. For avalanche generation, the stationary problem was treated in \cite{FN94,Mar85}, weak solutions on a short time interval in two dimensions were constructed in \cite{Nau95}, and local well-posedness in three dimensions was proved by maximal regularity in \cite{DR19}. The local theory of \cite{GHL26} rests on the compactness techniques of \cite{Isl03,Isl06,Nau95}.

The model of Degond, Lucquin-Desreux and Morrow \cite{DL07,Mor85}, analysed by Suzuki and
Tani \cite{ST18,ST20} and by Strauss and Suzuki \cite{SS19,SS22,SS25,SS26}, is the gas-discharge model
closest to \eqref{eq:sys}. It couples the Poisson equation to a transport equation for the
ions, which drift without diffusing, and to a drift--diffusion equation for the electrons, and
its source is that of \eqref{eq:F} with $\eta_0=0$: it is nonnegative, and the rate per
electron grows linearly in the field. The densities vanish on the electrodes, except that the
ions leave the cathode freely and that, in \cite{ST18} and in the stationary problems of \cite{SS22,SS26}, electrons are emitted at the cathode in proportion to the ion flux; every source of carriers vanishes in the absence of carriers, so the
non-ionized gas is a stationary solution, and breakdown is its loss of stability. In \cite{ST18} a unique
local-in-time solution is constructed in the half-space outside a compact cathode. Between
parallel plates, \cite{ST20} proves that the non-ionized state changes from asymptotically
stable, with global solutions for small initial data, to unstable at a sparking voltage, at
which non-trivial stationary solutions bifurcate from it, and \cite{SS19} continues this branch
until the densities become unbounded or the branch returns to the non-ionized state at a
second voltage. In annular regions between star-shaped surfaces in dimension two or three,
\cite{SS25} defines the sparking voltage through the lowest eigenvalue of a Schr\"odinger
operator, proves the same stability alternative for small data, and obtains the branch of
stationary solutions, globally in the radial case. Problem \eqref{eq:sys}--\eqref{eq:bc} differs from this model in four respects: the ions
diffuse; the net Townsend coefficient contains the attachment term, which these works omit; the
electrode data (A3) prescribe positive densities on both electrodes, so that carriers are
injected at the electrodes and a non-ionized state does not exist; and the boundary contains
the artificial Neumann part $\GN$, which meets the electrodes at corners.

\subsection{Contributions}\label{ssec:novelty}
We summarize the contributions of the paper.

\begin{enumerate}[leftmargin=2.4em,label=\textup{(\arabic*)}]
\item \emph{A continuation criterion in a weak norm.} The quadratic term of \eqref{eq:GHL-riccati} arises because the field in the ionization term is bounded through the $L^2$-norm of the densities (Remark~\ref{rem:origin}). Keeping the field in $L^2$ instead gives an inequality linear in the $L^2$ energy (Proposition~\ref{prop:field-controlled}), and hence Theorem~\ref{thm:A}. Unlike the usual continuation criteria for parabolic systems, stated in $L^\infty$- or $L^p$-norms of the unknowns, Theorem~\ref{thm:A} does not involve the total density and involves the charge only through its $\HDd$-norm, which is weaker than every $L^r$-norm with $r>2d/(d+2)$.

\item \emph{Global existence for an unbounded ionization rate.} The global results for the van Roosbroeck system and for the gas-discharge model of \cite{Isl03,Isl06} rest on generation rates per carrier that are bounded in the field. For the avalanche generation of semiconductor theory, whose source has the same quadratic structure as \eqref{eq:F}, the known results are stationary or local in time, and for the model of Degond, Lucquin-Desreux and Morrow the known global solutions start near the non-ionized state \cite{SS25,ST20} (Section~\ref{ssec:background}). Theorem~\ref{thm:B} requires no smallness of the voltage or of the data; to our knowledge, it is the first global existence result without smallness
assumptions for a drift--diffusion--Poisson model of gas discharge whose ionization rate is
unbounded in the field, and, as Remark~\ref{rem:constants} records, its estimates use no
property of $\Gamma$ beyond \eqref{eq:Gamma-bound}.

\item \emph{A closure without the free energy.} The proof of Proposition~\ref{prop:closure}, the bound on the electrostatic energy to which Theorem~\ref{thm:A} reduces global existence, uses neither the $L^2$ energy of \cite{GHL26} nor the free energy of the van Roosbroeck system (entropy plus electrostatic energy), whose dissipation does not control the Joule heating $\int_\Om(\mu_+p+\mu_-n)\abs{\nabla(\phi-\phi_D)}^2\dx$ (Remark~\ref{rem:freeenergy}). It uses the electrostatic energy of the space charge, $\ee:=\epsilon_0\norm{\nabla(\phi-\phi_D)}_{L^2}^2/2=\norm{p-n}_{\HDd}^2/(2\epsilon_0)$, and the $L^1$-distance of the densities from the electrode data, which controls the number of carriers; their sum satisfies a \emph{linear} Gronwall inequality, for two reasons. First, the ionization source cancels in the balance law of $\ee$, since carriers are created in pairs (Lemma~\ref{lem:energy-identity}). Second, the Joule heating is dissipative, because the densities are nonnegative, and, together with the number of carriers, it bounds both the battery term, through which the voltage source acts on $\ee$ (Corollary~\ref{cor:energy-ineq}), and, by \eqref{eq:source-Joule}, the ionization source of the $L^1$ estimate of Lemma~\ref{lem:L1}. Hypothesis (H1) enters only by making the ambipolar term of \eqref{eq:energy-identity} vanish, and (H2) only through these two bounds. The argument rests on the choice of these two quantities: with respect to them, the ionization
source, which produces the quadratic term in \eqref{eq:GHL-riccati}, either cancels or is dominated
by the Joule heating that the energy balance dissipates. Hence the Gronwall inequality is linear, and global existence does not depend on the size of the voltage or of the data.

\end{enumerate}

\subsection{Strategy and outline}\label{ssec:strategy}
The proofs use three properties of \eqref{eq:sys}. First, the ionization source $F=\Gamma(|\nabla\phi|)\,n$ is linear in $n$, with a rate that satisfies $\abs{\Gamma(|\nabla\phi|)}\le\mu_-\gamma_*|\nabla\phi|$ by \eqref{eq:Gamma-bound}. Second, the same source $F$ enters \eqref{eq:p} and \eqref{eq:n}, so that ionization creates carriers, not charge \cite[Remark~1.3]{GHL26}: subtracting the two equations and using $\nabla\cdot\bv=0$ gives, for the charge $q:=p-n$ and the conductivity $g:=\mu_+p+\mu_-n$,
\begin{equation}\label{eq:charge-eq}
\partial_tq+\bv\cdot\nabla q=\nabla\cdot\big(\epsilon_+\nabla p-\epsilon_-\nabla n+g\,\nabla\phi\big).
\end{equation}
Third, $p,n\ge0$ \cite[Lemma~4.1]{GHL26}.

Section~\ref{sec:setting} collects notation, assumptions, the notion of weak solution and some auxiliary inequalities. Section~\ref{sec:contcrit} proves Theorem~\ref{thm:A}. In the $L^2$ energy identity (Lemma~\ref{lem:L2identity}) the drift terms contribute only data terms and the screening term $\frac{\mu_+\mu_-}{2\epsilon_0}\int_\Om(p+n)(p-n)^2\dx$, which is nonnegative by the third property (as already used, for the truncated problem, in \cite[proof of Theorem~4.2]{GHL26}). By the first property and interpolation (Lemma~\ref{lem:source}(ii)), the $L^2$-norms of the densities then grow at most like $\exp\big(C\int_0^t(1+\norm{\nabla\phi}_{L^2}^{2\kappa})\dr\big)$ (Proposition~\ref{prop:field-controlled}), and the blow-up alternative of Lemma~\ref{lem:alternative} concludes. Since $\norm{\nabla\phi}_{L^2}\le\norm{\nabla\phi_D}_{L^2}+(2\ee/\epsilon_0)^{1/2}$, Theorem~\ref{thm:A} reduces global existence to a bound on $\ee$ over bounded time intervals. Section~\ref{sec:global} obtains such a bound under (H1)--(H2), and with it Theorem~\ref{thm:B}: testing \eqref{eq:charge-eq} with $\phi-\phi_D$ gives the balance law of Lemma~\ref{lem:energy-identity} for $\ee$, and the energy inequality of Corollary~\ref{cor:energy-ineq} closes against the $L^1$ estimate of Lemma~\ref{lem:L1}, which controls the number of carriers (Proposition~\ref{prop:closure}). Remark~\ref{rem:H1} and Proposition~\ref{prop:conditional} record what the argument gives without (H1). Section~\ref{sec:numerics} reports two numerical experiments on a discharge in a micrometre gap, with the Townsend constants and the mobilities of air: the dependence on the voltage, in connection with the breakdown threshold, and the terms of \eqref{eq:needed} in a simulation with the diffusivities of air. Section~\ref{sec:remarks} discusses unequal diffusivities, the failure of (H2) at touch-down, the singularities that Theorems~\ref{thm:A} and~\ref{thm:B} leave possible, and the breakdown threshold. Appendix~\ref{app:scheme} describes the numerical scheme.

\section{Preliminaries}\label{sec:setting}

\subsection{Notation}\label{ssec:notation}
For $1\le r\le\infty$ we write $\norm{\cdot}_r:=\norm{\cdot}_{L^r(\Om)}$, and $\norm{\cdot}_{H^k}:=\norm{\cdot}_{H^k(\Om)}$. We set
\[
\HD:=\{u\in H^1(\Om):\ u=0\ \text{on }\GD\},\qquad \norm{u}_{\HD}:=\norm{\nabla u}_2,
\]
which is a norm by the Poincar\'e inequality $\norm{u}_2\le C_P\norm{\nabla u}_2$ on $\HD$, valid since $\GD$ has positive surface measure by (A1) below, and denote by $\HDd$ the dual of $\HD$ with the pairing $\pair{\cdot,\cdot}$; $L^2(\Om)$ is identified with a subspace of $\HDd$. For $u\in L^2((0,T),\HD)\cap H^1((0,T),\HDd)$ we recall that $u\in C([0,T],L^2(\Om))$ and $\frac{d}{dt}\norm{u}_2^2=2\pair{\partial_tu,u}$ a.e.\ in $(0,T)$ \cite[\S5.9]{Eva10}.

Throughout the paper we use the abbreviations
\begin{equation}\label{eq:abbrev}
\begin{gathered}
q:=p-n,\qquad s:=p+n,\qquad g:=\mu_+p+\mu_-n,\\
\tp:=p-p_D,\qquad \tn:=n-n_D .
\end{gathered}
\end{equation}
For $p,n\ge0$, these satisfy $\abs{q}\le s$ and $\min\{\mu_+,\mu_-\}\,s\le g\le\max\{\mu_+,\mu_-\}\,s$.

\subsection{Assumptions}
The following is the standing assumption \cite[Assumption~1.5]{GHL26}, reorganized, with the velocity field prescribed and bounded on the whole half-line; it also states two conditions that belong to the setting of \cite[\S1]{GHL26}: the tangency
$\bnu\cdot\bv=0$ on $\partial\Om$ and the positive surface measure of $\GD$.

\begin{assumption}\label{ass:main}\hfill
\begin{itemize}[leftmargin=2.6em]
\item[\textup{(A1)}] \textup{(Domain.)} $\Om\subset\R^d$, $d\in\{2,3\}$, is a bounded domain, $\partial\Om=\overline{\GD}\cup\overline{\GN}$ with $\GD\cap\GN=\emptyset$, $\GD$ and $\GN$ are open subsets of $\partial\Om$ of class $C^{1,1}$ enclosing an angle at most $\pi/2$, and $\GD$ has positive surface measure. In particular the elliptic regularity of Lemma~\ref{lem:aux}\textup{(ii)} below holds.
\item[\textup{(A2)}] \textup{(Electrode potential.)} The electrode potential is the trace on $\GD$ of a function in $H^2(\Om)$. We denote by $\phi_D\in H^2(\Om)$ its \emph{harmonic extension}, i.e.\ the unique solution of
\[
\Delta\phi_D=0\ \text{ in }\Om,\qquad \phi_D=\text{electrode potential on }\GD,\qquad \bnu\cdot\nabla\phi_D=0\ \text{ on }\GN,
\]
which belongs to $H^2(\Om)$ by \cite[Proposition~2.2]{GHL26} and satisfies $\int_\Om\nabla\phi_D\cdot\nabla\chi\dx=0$ for all $\chi\in\HD$.
\item[\textup{(A3)}] \textup{(Carrier data.)} $p_D,n_D\in W^{2,\sigma}(\Om)$ for some $\sigma\ge3$, $p_D,n_D\ge0$ in $\Om$, $\bnu\cdot\nabla p_D=\bnu\cdot\nabla n_D=0$ on $\GN$, and $R_a\le p_D,n_D\le R_b$ on $\GD$ for some constants $0<R_a\le R_b$. The initial data satisfy $p_0,n_0\in W^{2,\sigma}(\Om)$, $p_0,n_0\ge0$, $p_0=p_D$ and $n_0=n_D$ on $\GD$, $\bnu\cdot\nabla p_0=\bnu\cdot\nabla n_0=0$ on $\GN$.
\item[\textup{(A4)}] \textup{(Gas velocity.)} $\bv\in L^\infty\big((0,\infty),H^2(\Om)^d\big)$, $\nabla\cdot \bv=0$ in $\Om$ and $\bnu\cdot \bv=0$ on $\partial\Om$ for a.e.\ $t$.
\end{itemize}
\end{assumption}

By (A3), and since $W^{2,\sigma}(\Om)\hookrightarrow C^{0,\alpha}(\overline\Om)$ for $\sigma\ge3$ and $d\le3$, the data $p_D,n_D$ are nonnegative and bounded, $\nabla p_D,\nabla n_D\in W^{1,\sigma}(\Om)^d\hookrightarrow L^r(\Om)^d$ for every $r<\infty$, so that $\bnu\cdot\nabla p_D$ and $\bnu\cdot\nabla n_D$ have traces on $\partial\Om$, and $\Delta p_D,\Delta n_D\in L^\sigma(\Om)\subset L^1(\Om)$; together with the Neumann compatibility on $\GN$, these are the only properties of $p_D,n_D$ used in the estimates of Sections~\ref{sec:L2} and~\ref{sec:global}. Likewise $H^2(\Om)\hookrightarrow L^\infty(\Om)$, so $\bv\in L^\infty(\Om\times(0,\infty))$, and we write
\begin{equation}\label{eq:v-infty}
\norm{\bv}_\infty:=\norm{\bv}_{L^\infty(\Om\times(0,\infty))}\le C\norm{\bv}_{L^\infty((0,\infty),H^2)},
\end{equation}
with $C$ the constant of that embedding, which depends only on $\Om$. The velocity enters the estimates of Sections~\ref{sec:L2} and~\ref{sec:global} only through this number. We write
\begin{equation}\label{eq:data}
\begin{gathered}
\mathcal D_0:=\big(d,\Om,\epsilon_\pm,\mu_\pm,\epsilon_0,\alpha_1,\alpha_2,\eta_0,\norm{\phi_D}_{H^2},\norm{p_D}_{W^{2,\sigma}},\norm{n_D}_{W^{2,\sigma}}\big),\\
\mathcal D:=\big(\mathcal D_0,\norm{\bv}_{L^\infty((0,\infty),H^2)}\big),
\end{gathered}
\end{equation}
and, when (H2) is in force, we include $\norm{\nabla\phi_D}_\infty$ in $\mathcal D$. The letter $C$ denotes a positive constant depending only on $\mathcal D$; it may change from line to line, but it never depends on $T$ or $t$. Every dependence on time is written out explicitly.

\subsection{Weak solutions}\label{ssec:weak}
The following definition is equivalent to \cite[Definition~1.6]{GHL26}; we use the conservative form of the parabolic equations, which is obtained from the form in \cite{GHL26} by an integration by parts that uses $\phi\in H^2(\Om)$, $\bnu\cdot\nabla\phi=0$ on $\GN$, the vanishing of the test functions  on $\GD$, $\nabla\cdot \bv=0$ and $\bnu\cdot \bv=0$.

\begin{definition}\label{def:weak}
Let $T>0$. A triple $(\phi,p,n)$ is a \emph{weak solution} of \eqref{eq:sys}--\eqref{eq:bc} on $(0,T)$ if
\begin{enumerate}[leftmargin=2.2em,label=\textup{(\roman*)}]
\item $\phi\in L^\infty((0,T),H^2(\Om))$, $p,n\in L^\infty((0,T),L^2(\Om))\cap L^2((0,T),H^1(\Om))\cap H^1((0,T),\HDd)$, $p-p_D,\ n-n_D\in L^2((0,T),\HD)$, and $p,n\ge0$ a.e.\ in $\Om\times(0,T)$;
\item for a.e.\ $t\in(0,T)$, $\psi(t):=\phi(t)-\phi_D\in\HD$ and
\begin{equation}\label{eq:weak-poisson}
\epsilon_0\int_\Om\nabla\psi(t)\cdot\nabla\chi\dx=\int_\Om q(t)\,\chi\dx\qquad\text{for all }\chi\in\HD;
\end{equation}
\item for a.e.\ $t\in(0,T)$ and all $f\in\HD$,
\begin{align}
\pair{\partial_tp(t),f}+\int_\Om\big(\epsilon_+\nabla p+\mu_+p\nabla\phi-p\bv\big)\cdot\nabla f\dx&=\int_\Om F\big(n,|\nabla\phi|\big)f\dx,\label{eq:weak-p}\\
\pair{\partial_tn(t),f}+\int_\Om\big(\epsilon_-\nabla n-\mu_-n\nabla\phi-n\bv\big)\cdot\nabla f\dx&=\int_\Om F\big(n,|\nabla\phi|\big)f\dx;\label{eq:weak-n}
\end{align}
\item $p(0)=p_0$ and $n(0)=n_0$ in $L^2(\Om)$.
\end{enumerate}
\end{definition}

All integrals in \eqref{eq:weak-p}--\eqref{eq:weak-n} are finite for a.e.\ $t$: by (A2) and Lemma~\ref{lem:aux}(ii), $\nabla\phi(t)\in L^6(\Om)$, while $p(t),n(t)\in H^1(\Om)\hookrightarrow L^6(\Om)$, so that $p\nabla\phi$, $F$ and $n\bv$ belong to $L^3(\Om)$. By \eqref{eq:weak-poisson} and (A2), $\phi(t)$ is the solution of the mixed problem for \eqref{eq:poisson} with datum $q(t)\in L^2(\Om)$, hence the $H^2$-regularity in (i) is automatic and $\bnu\cdot\nabla\phi=0$ on $\GN$ in the sense of traces. Moreover, since $q\mapsto\psi$ is a bounded linear map $\HDd\to\HD$ and $q\in H^1((0,T),\HDd)$,
\begin{equation}\label{eq:psi-regularity}
\begin{split}\psi\in H^1((0,T),\HD)\quad\text{and}\qquad\qquad\qquad\quad\quad\\ \epsilon_0\int_\Om\nabla\partial_t\psi(t)\cdot\nabla\chi\dx=\pair{\partial_tq(t),\chi}\quad\text{for all }\chi\in\HD\text{ and a.e.\ }t.\end{split}
\end{equation}

\subsection{Auxiliary inequalities}

\begin{lemma}\label{lem:aux}
Let \textup{(A1)} hold.
\begin{enumerate}[leftmargin=2.2em,label=\textup{(\roman*)}]
\item \textup{(Gagliardo--Nirenberg--Ladyzhenskaya.)} There is $C_{GN}$ such that for all $u\in H^1(\Om)$
\[
\norm{u}_4^2\le C_{GN}\norm{u}_2\norm{u}_{H^1}\ \ (d=2),\qquad
\norm{u}_4^2\le C_{GN}\norm{u}_2^{1/2}\norm{u}_{H^1}^{3/2}\ \ (d=3),
\]
and $\norm{u}_3+\norm{u}_6\le C_{GN}\norm{u}_{H^1}$. Moreover $H^2(\Om)\hookrightarrow W^{1,6}(\Om)\cap L^\infty(\Om)$.
\item \textup{(Poisson problem.)} For $f\in\HDd$ there is a unique $u\in\HD$ with $\epsilon_0\int_\Om\nabla u\cdot\nabla\chi\dx=\pair{f,\chi}$ for all $\chi\in\HD$, and $\norm{\nabla u}_2=\epsilon_0^{-1}\norm{f}_{\HDd}$. If $f\in L^2(\Om)$ then $u\in H^2(\Om)$, $\bnu\cdot\nabla u=0$ on $\GN$, and $\norm{u}_{H^2}\le C_{\mathrm{ell}}\norm{f}_2$.
\item \textup{(Chain rule.)} Let $u\in L^2((0,T),\HD)\cap H^1((0,T),\HDd)$, and let $\Sigma\in C^1(\R)$ be convex with $\Sigma(0)=\Sigma'(0)=0$ and $\Sigma'$ Lipschitz continuous. Then $t\mapsto\int_\Om\Sigma(u(t))\dx$ is absolutely continuous on $[0,T]$ and
\[\frac{d}{dt}\int_\Om\Sigma(u(t))\dx=\pair{\partial_tu(t),\Sigma'(u(t))}\qquad\text{for a.e.\ }t\in(0,T).\]
\end{enumerate}
\end{lemma}

\begin{proof}
(i) is the Gagliardo--Nirenberg inequality \cite{Nir59} on $\R^d$ combined with a bounded extension operator $H^1(\Om)\to H^1(\R^d)$, together with the Sobolev embeddings for $d\le3$.

(ii) is the Lax--Milgram lemma; the identity holds because, by the Cauchy--Schwarz inequality,
\[
\norm{f}_{\HDd}=\sup_{\substack{\chi\in\HD,\\ \norm{\nabla\chi}_2=1}}\epsilon_0\int_\Om\nabla u\cdot\nabla\chi\dx=\epsilon_0\norm{\nabla u}_2 ,
\]
the supremum being attained at $\chi=u/\norm{\nabla u}_2$ if $u\neq0$. The regularity statement is \cite[Proposition~2.2]{GHL26}, and the Neumann condition follows from the weak formulation once $u\in H^2(\Om)$.

(iii) is \cite[Lemma~1.5]{AL83}; see also \cite[Chapter~7]{Rou13}. Moreover $0\le\Sigma(u)\le\tfrac L2u^2$, with $L$ the Lipschitz constant of $\Sigma'$, and $\Sigma'(u(t))\in\HD$ because $\Sigma'$ is Lipschitz with $\Sigma'(0)=0$.
\end{proof}

\section{Continuation criterion: proof of Theorem~\ref{thm:A}}\label{sec:contcrit}

\subsection{The local theory and the maximal solution}\label{ssec:local}
The local theory is \cite[Theorem~1.1]{GHL26}; the system, the source \eqref{eq:F} and the notion of weak solution are those of \cite{GHL26} (Section~\ref{ssec:weak}). The proofs below use \cite{GHL26} only through Theorem~\ref{thm:GHL}, Proposition~\ref{prop:restart}, the positivity $p,n\ge0$ in Definition~\ref{def:weak}(i) and the regularity statement of Lemma~\ref{lem:aux}(ii); the remaining references to \cite{GHL26} serve to locate the origin of the quadratic term in \eqref{eq:GHL-riccati} and to compare hypotheses. We use this theorem in the following form.

\begin{theorem}[{cf.\ \cite[Theorem~1.1]{GHL26}}]\label{thm:GHL}
Let \textup{(A1)}, \textup{(A2)} and the conditions on $p_D,n_D$ in \textup{(A3)} hold. For all $R,\Lambda\ge0$ there is $\tau_{\mathcal D_0}(R,\Lambda)>0$, depending only on $\mathcal D_0$, $R$ and $\Lambda$, with the following property: for all initial data $p_0,n_0$ satisfying \textup{(A3)} with $\norm{p_0}_2+\norm{n_0}_2\le R$ and every velocity field $\bv$ satisfying \textup{(A4)} with $\norm{\bv}_{L^\infty((0,\infty),H^2)}\le\Lambda$, problem \eqref{eq:sys}--\eqref{eq:bc} has a weak solution on $\big(0,\tau_{\mathcal D_0}(R,\Lambda)\big)$. Moreover, two weak solutions on a common interval $(0,T')$ with the same $p_0,n_0$ and $\bv$ coincide on $(0,T')$.
\end{theorem} 

The existence time is a time $T_1$ up to which the denominator of the bound \cite[(4.15)]{GHL26} stays positive. This denominator is nonincreasing in the constants $\tilde H_4,\tilde H_5,\tilde H_6$ in the proof of \cite[Theorem~4.2]{GHL26}, and these depend on the initial data only through $\norm{p_0}_2^2+\norm{n_0}_2^2$ and on the velocity only through $\esssup_t\norm{\bv(t)}_{H^2}$, as nondecreasing functions of these two quantities with coefficients depending only on $\mathcal D_0$. Hence $T_1$ can be chosen, for all initial data with $\norm{p_0}_2+\norm{n_0}_2\le R$ and all velocity fields with $\norm{\bv}_{L^\infty((0,\infty),H^2)}\le\Lambda$, as one positive number depending only on $\mathcal D_0$, $R$ and $\Lambda$; this number is the time $\tau_{\mathcal D_0}(R,\Lambda)$ of Theorem~\ref{thm:GHL}. This uniformity is what allows the problem to be restarted: in the proof of Lemma~\ref{lem:alternative} below, the solution is continued beyond a time $t_0$ by solving the problem with the data $(p(t_0),n(t_0))$ and the velocity field $\bv(\cdot+t_0)$, which satisfies (A4) with the same $\Lambda$. By Definition~\ref{def:weak}(i), however, the data $(p(t_0),n(t_0))$ have only the regularity \eqref{eq:restart-data} below, which is less than (A3) requires. The following proposition removes this obstacle.


\begin{proposition}[Restart]\label{prop:restart}
Theorem~\ref{thm:GHL}, with the same existence time $\tau:=\tau_{\mathcal D_0}(R,\Lambda)$, remains valid if the conditions on $(p_0,n_0)$ in \textup{(A3)} are replaced by
\begin{equation}\label{eq:restart-data}
p_0,n_0\in H^1(\Om),\qquad p_0,n_0\ge0\ \text{a.e.\ in }\Om,\qquad p_0-p_D,\ n_0-n_D\in\HD .
\end{equation}
\end{proposition}

\begin{proof}
In \cite{GHL26} the initial data enter the a priori bounds of \cite[Theorem~4.2]{GHL26}, and
hence the existence time, only through $\norm{p_0}_2^2+\norm{n_0}_2^2$, and the regularity required of $(p_0,n_0)$ in (A3) is used only in the estimates of \cite[Corollary~4.3 and Remark~4.4]{GHL26} and in the corresponding estimate in \cite[Appendix~A]{GHL26}. These estimates serve to pass to the limit in the nonlinear terms and in the truncated Poisson equation, as $M\to\infty$ in \cite[\S5.1]{GHL26} and as $j\to\infty$ in \cite[Appendix~A]{GHL26}. For this, the a priori $L^2$ bounds suffice, together with the strong $L^2$-convergence of the approximations, which the Aubin--Lions lemma gives from these bounds and the resulting bound on the time derivatives. Indeed, the cut-off of the charge in the truncated Poisson equation is $1$-Lipschitz and the identity on $[-M,M]$, so that, as $M\to\infty$, it converges in $L^2(\Om\times(0,\tau))$ together with its argument, and Lemma~\ref{lem:aux}(ii) then gives the convergence of the potentials in $L^2((0,\tau),H^1(\Om))$; by interpolation with the bounds on the densities in $L^2((0,\tau),H^1(\Om))$ and on the potentials in $L^\infty((0,\tau),H^2(\Om))$, the densities converge in $L^2((0,\tau),L^3(\Om))$ and the gradients of the potentials in $L^4((0,\tau),L^3(\Om))$, which suffices for the limit of the products $p\nabla\phi$, $n\nabla\phi$ and $\Gamma(|\nabla\phi|)\,n$. The uniqueness part is unaffected: the proof in \cite[\S5.2]{GHL26} compares two weak solutions in the class of Definition~\ref{def:weak} and uses no property of the initial data beyond $p_0,n_0\in L^2(\Om)$.
\end{proof}

Let Assumption~\ref{ass:main} hold, and let $\mathcal T$ be the set of $T'>0$ such that \eqref{eq:sys}--\eqref{eq:bc} has a weak solution on $(0,T')$. By Theorem~\ref{thm:GHL}, $\mathcal T\neq\emptyset$; by uniqueness, weak solutions on $(0,T')$ and on $(0,T'')$ coincide on $(0,\min\{T',T''\})$. Hence there is a unique triple $(\phi,p,n)$ on $(0,T_{\max})$, $T_{\max}:=\sup\mathcal T\in(0,\infty]$, whose restriction to $(0,T')$ is a weak solution for every $T'<T_{\max}$. We call it the \emph{maximal weak solution}. By Proposition~\ref{prop:restart}, the same construction applies to initial data that satisfy \eqref{eq:restart-data} in place of the conditions of (A3), and Lemma~\ref{lem:alternative} and Theorems~\ref{thm:A} and~\ref{thm:B} hold for such data as well, since their proofs use the initial data only through Proposition~\ref{prop:restart} and through the $L^2$-norms of $p_0,n_0$, which are finite for data satisfying \eqref{eq:restart-data}. The following blow-up alternative is the continuation property conjectured in \cite[Remark~1.2]{GHL26}.

\begin{lemma}[Blow-up alternative]\label{lem:alternative}
Let Assumption~\ref{ass:main} hold and let $(\phi,p,n)$ be the maximal weak solution. If $T_{\max}<\infty$, then
\[
\limsup_{t\rightarrow T_{\max}}\big(\norm{p(t)}_2+\norm{n(t)}_2\big)=\infty .
\]
\end{lemma}

\begin{figure}[t]
\centering
\begin{tikzpicture}[x=1.08cm,y=0.78cm,font=\small]
  \def\tB{4} \def\tC{5.5} \def\tD{7.5} \def\tE{9} \def\tF{10.5}
  \fill[black!10] (\tC,0) rectangle (\tD,3.45);
  \foreach \x in {\tC,\tD,\tE,\tF} \draw[densely dotted,black!45] (\x,0) -- (\x,3.45);
  \draw[line width=2.4pt] (0,2.9) -- (\tD,2.9);
  \draw[line width=1pt,dashed] (\tD,2.9) -- (\tE,2.9);
  \node[anchor=south west,inner sep=1pt] at (0,3.05) {$(\phi,p,n)$};
  \draw[line width=2.4pt,black!50] (\tC,1.85) -- (\tF,1.85);
  \node[anchor=east] at (\tC-0.1,1.85) {$(\hat\phi,\hat p,\hat n)(\cdot-t_0)$};
  \draw[|<->|,>=stealth] (\tC,2.25) -- (\tF,2.25) node[midway,above=1pt,inner sep=1pt] {$\tau$};
  \draw[line width=2.4pt] (0,0.8) -- (\tD,0.8);
  \draw[line width=2.4pt,black!50] (\tD,0.8) -- (\tF,0.8);
  \node[anchor=south west,inner sep=1pt] at (0,0.95) {$(\bar\phi,\bar p,\bar n)$};
  \draw[->] (-0.2,0) -- (\tF+0.9,0) node[right] {$t$};
  \foreach \x/\lab in {0/0,\tB/T_{\max}-\tau,\tC/t_0,\tD/T',\tE/T_{\max},\tF/t_0+\tau}
    \draw (\x,0.1) -- (\x,-0.1) node[below] {$\lab$};
\end{tikzpicture}
\caption{The restart in the proof of Lemma~\ref{lem:alternative}. The maximal weak solution is used on $(0,T')$ and the restarted solution, shifted by $t_0$, on $(t_0,t_0+\tau)$; the two coincide on the overlap $(t_0,T')$ (shaded). The glued triple $(\bar\phi,\bar p,\bar n)$ is a weak solution on $(0,t_0+\tau)$, and $t_0+\tau>T_{\max}$.}
\label{fig:restart}
\end{figure}
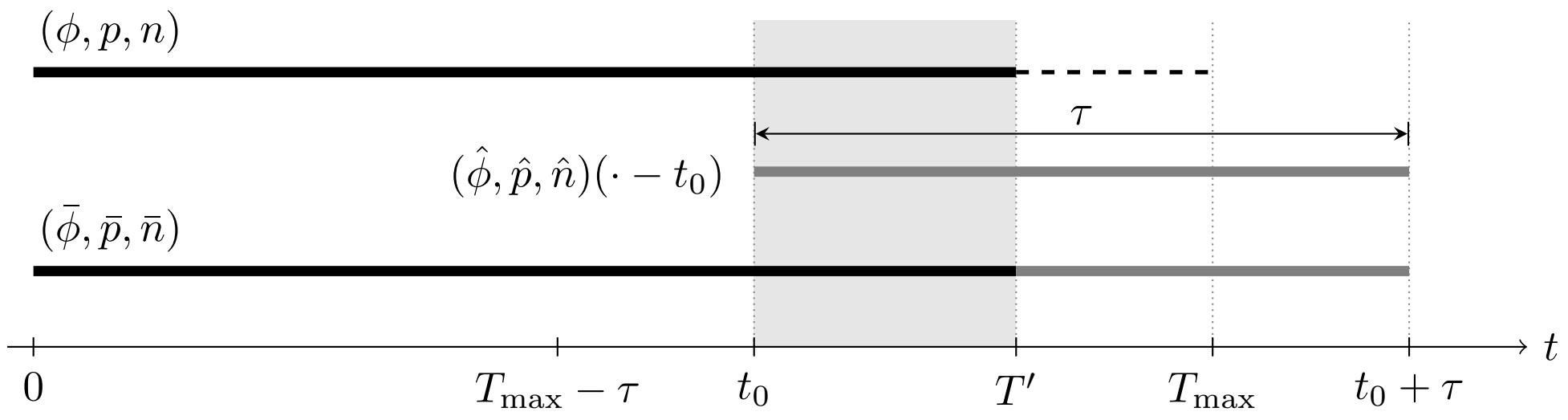

\begin{proof}
Assume that $T_{\max}<\infty$ and  \begin{equation}\label{finite_sup_L2}\limsup_{t\to T_{\max}}\big(\norm{p(t)}_2+\norm{n(t)}_2\big)<\infty.\end{equation}
For every $S\in(0,T_{\max})$, the maximal weak solution is a weak solution on $(0,S)$, so that $p-p_D,\ n-n_D\in C([0,S],L^2(\Om))$ (Section~\ref{ssec:notation}). This, together with \eqref{finite_sup_L2}, gives
\[
R:=\sup_{t\in(0, T_{\max})}\big(\norm{p(t)}_2+\norm{n(t)}_2\big)<\infty .
\]
Put $\Lambda:=\norm{\bv}_{L^\infty((0,\infty),H^2)}$ and $\tau:=\tau_{\mathcal D_0}(R,\Lambda)$. By Definition~\ref{def:weak}(i), for a.e.\ $t\in(0,T_{\max})$ the data $(p(t),n(t))$ satisfy \eqref{eq:restart-data}, that is, $p(t),n(t)\in H^1(\Om)$, $p(t),n(t)\ge0$ a.e.\ in $\Om$ and $p(t)-p_D,\ n(t)-n_D\in\HD$. We fix such a time $t_0\in(\max\{0,T_{\max}-\tau\},T_{\max})$ and then $T'\in(t_0,T_{\max})$, so that $T'-t_0<\tau$ and $t_0+\tau>T_{\max}$ (Figure~\ref{fig:restart}). Since $\bv(\cdot+t_0)$ satisfies (A4) with $\norm{\bv(\cdot+t_0)}_{L^\infty((0,\infty),H^2)}\le\Lambda$ and $\norm{p(t_0)}_2+\norm{n(t_0)}_2\le R$, Proposition~\ref{prop:restart} gives a weak solution $(\hat\phi,\hat p,\hat n)$ on $(0,\tau)$ of \eqref{eq:sys}--\eqref{eq:bc} with the velocity field $\bv(\cdot+t_0)$ and the initial data $(p(t_0),n(t_0))$. The shifted triple $(\phi,p,n)(\cdot+t_0)$ is a weak solution of the same problem on $(0,T'-t_0)\subset(0,\tau)$, so that, by the uniqueness part of Proposition~\ref{prop:restart},
\[
(\hat\phi,\hat p,\hat n)=(\phi,p,n)(\cdot+t_0)\qquad\text{on }(0,T'-t_0).
\]
Hence the triple
\[
(\bar\phi,\bar p,\bar n)(t):=
\begin{cases}
(\phi,p,n)(t), & t\in(0,T'),\\
(\hat\phi,\hat p,\hat n)(t-t_0), & t\in(t_0,t_0+\tau),
\end{cases}
\]
is well defined, and it is a weak solution on $(0,t_0+\tau)$. Indeed, $(\phi,p,n)$ satisfies conditions (ii) and (iii) of Definition~\ref{def:weak} for a.e.\ $t\in(0,T')$, and so does $(\hat\phi,\hat p,\hat n)(\cdot-t_0)$ for a.e.\ $t\in(t_0,t_0+\tau)$, because $(\hat\phi,\hat p,\hat n)$ solves the problem with the velocity field $\bv(\cdot+t_0)$; hence $(\bar\phi,\bar p,\bar n)$ satisfies (ii) and (iii) for a.e.\ $t\in(0,t_0+\tau)$. The integrability and sign conditions in Definition~\ref{def:weak}(i) for $(\bar\phi,\bar p,\bar n)$ on $(0,t_0+\tau)$ follow from those for $(\phi,p,n)$ on $(0,T')$ and for $(\hat\phi,\hat p,\hat n)(\cdot-t_0)$ on $(t_0,t_0+\tau)$. On the overlap $(t_0,T')$, the weak time derivatives $\partial_tp$ and $(\partial_t\hat p)(\cdot-t_0)$ coincide, and so do $\partial_tn$ and $(\partial_t\hat n)(\cdot-t_0)$. By a partition of unity in time, the weak time derivatives on the two intervals therefore combine into weak time derivatives of $\bar p$ and $\bar n$ on $(0,t_0+\tau)$, so that $\bar p,\bar n\in H^1((0,t_0+\tau),\HDd)$. Finally, (iv) holds since $(\bar\phi,\bar p,\bar n)$ coincides with $(\phi,p,n)$ on $(0,T')$. As $t_0+\tau>T_{\max}$, the existence of this weak solution contradicts the definition of $T_{\max}$.
\end{proof}

\subsection{A field-controlled \texorpdfstring{$L^2$}{L2} estimate}\label{sec:L2}

Throughout this subsection $(\phi,p,n)$ is a weak solution on $(0,T)$ under Assumption~\ref{ass:main}, with the abbreviations \eqref{eq:abbrev}. We set
\begin{equation}\label{eq:Y}
\Y(t):=\mu_-\norm{\tp(t)}_2^2+\mu_+\norm{\tn(t)}_2^2,\qquad
\theta:=\min\{\epsilon_+\mu_-,\ \epsilon_-\mu_+\},
\end{equation}
so that $\Y$ is a weighted analogue, for the differences $\tp,\tn$, of $Y$ in \eqref{eq:GHL-riccati}. Since $\abs{\norm{p}_2-\norm{\tp}_2}\le\norm{p_D}_2$ and $\abs{\norm{n}_2-\norm{\tn}_2}\le\norm{n_D}_2$, we have $\norm{p}_2^2+\norm{n}_2^2\le C(1+\Y)$ and $\Y\le C(1+\norm{p}_2^2+\norm{n}_2^2)$.

\begin{lemma}[$L^2$ energy identity]\label{lem:L2identity}
For a.e.\ $t\in(0,T)$,
\begin{equation}\label{eq:L2identity}
\frac12\frac{d\Y}{dt}+\epsilon_+\mu_-\norm{\nabla\tp}_2^2+\epsilon_-\mu_+\norm{\nabla\tn}_2^2+\frac{\mu_+\mu_-}{2\epsilon_0}\int_\Om s\,q^2\dx=\Rdata+\RF ,
\end{equation}
where
\begin{align*}
\RF&:=\int_\Om F\big(n,|\nabla\phi|\big)\big(\mu_-\tp+\mu_+\tn\big)\dx,\\
\Rdata&:=-\epsilon_+\mu_-\int_\Om\nabla p_D\cdot\nabla\tp\dx-\epsilon_-\mu_+\int_\Om\nabla n_D\cdot\nabla\tn\dx
-\mu_+\mu_-\int_\Om\big(p_D\nabla\tp-n_D\nabla\tn\big)\cdot\nabla\phi\dx\\
&\qquad+\mu_-\int_\Om p_D\,\bv\cdot\nabla\tp\dx+\mu_+\int_\Om n_D\,\bv\cdot\nabla\tn\dx\\
&\qquad+\frac{\mu_+\mu_-}{\epsilon_0}\int_\Om\big(pp_D-nn_D\big)q\dx-\frac{\mu_+\mu_-}{2\epsilon_0}\int_\Om\big(p_D^2-n_D^2\big)q\dx .
\end{align*}
\end{lemma}

\begin{proof}
The identity follows by testing the two parabolic equations \eqref{eq:weak-p}, \eqref{eq:weak-n} and integrating by parts. These integrations by parts are justified because, for a.e.\ $t$, $\tp^2(t)$, $\tn^2(t)\in W^{1,3/2}(\Om)$, $\phi(t)\in H^2(\Om)$ and $\bv(t)\in H^2(\Om)^d$. They produce the boundary terms $\frac12\int_{\partial\Om}(\bnu\cdot\bv)\,\tp^2\dS$ and $\frac12\int_{\partial\Om}(\bnu\cdot\nabla\phi)\,\tp^2\dS$, and the same terms with $\tn$; these vanish because $\bnu\cdot\bv=0$ on $\partial\Om$, $\tp=\tn=0$ on $\GD$ and $\bnu\cdot\nabla\phi=0$ on $\GN$. The crucial step  is the algebra in the drift terms, which produces the screening term $\frac{\mu_+\mu_-}{2\epsilon_0}\int_\Om s\,q^2\dx$ on the left of \eqref{eq:L2identity}. We record the whole computation.

We take $f=\mu_-\tp$ in \eqref{eq:weak-p} and $f=\mu_+\tn$ in \eqref{eq:weak-n}, write $p=\tp+p_D$ and $n=\tn+n_D$ in the fluxes, and add the two equations. By Lemma~\ref{lem:aux}(iii), applied with $\Sigma(u)=u^2/2$ to $\tp$ and $\tn$, $\Y$ is absolutely continuous on $[0,T]$; since $p_D,n_D$ do not depend on $t$, we have $\partial_tp=\partial_t\tp$ and $\partial_tn=\partial_t\tn$, so that the time derivatives contribute $\frac12\frac{d\Y}{dt}$. The source terms give $\RF$ on the right, and the remaining terms on the left are as follows.

\emph{Diffusion.} We have
\begin{multline*}
\epsilon_+\mu_-\int_\Om\nabla p\cdot\nabla\tp\dx+\epsilon_-\mu_+\int_\Om\nabla n\cdot\nabla\tn\dx\\
=\epsilon_+\mu_-\norm{\nabla\tp}_2^2+\epsilon_-\mu_+\norm{\nabla\tn}_2^2
+\epsilon_+\mu_-\int_\Om\nabla p_D\cdot\nabla\tp\dx+\epsilon_-\mu_+\int_\Om\nabla n_D\cdot\nabla\tn\dx .
\end{multline*}

\emph{Convection.} We write $p=\tp+p_D$ and $n=\tn+n_D$. Since $\nabla\cdot\bv=0$ in $\Om$ and $\bnu\cdot\bv=0$ on $\partial\Om$, an integration by parts gives
\[
\int_\Om\tp\,\bv\cdot\nabla\tp\dx=\frac12\int_\Om\bv\cdot\nabla(\tp^2)\dx=0,
\]
and likewise for $\tn$. Hence only the terms with $p_D$ and $n_D$ remain:
\[
-\mu_-\int_\Om p\,\bv\cdot\nabla\tp\dx-\mu_+\int_\Om n\,\bv\cdot\nabla\tn\dx
=-\mu_-\int_\Om p_D\,\bv\cdot\nabla\tp\dx-\mu_+\int_\Om n_D\,\bv\cdot\nabla\tn\dx .
\]

\emph{Drift.} By Green's formula and \eqref{eq:poisson},
\[
\int_\Om\tp\,\nabla\phi\cdot\nabla\tp\dx=\frac1{2\epsilon_0}\int_\Om\tp^2q\dx,
\]
and likewise for $\tn$, whose drift enters with the opposite sign. Hence
\[
\mu_+\mu_-\int_\Om\big(p\,\nabla\tp-n\,\nabla\tn\big)\cdot\nabla\phi\dx
=\frac{\mu_+\mu_-}{2\epsilon_0}\int_\Om\big(\tp^2-\tn^2\big)q\dx
+\mu_+\mu_-\int_\Om\big(p_D\nabla\tp-n_D\nabla\tn\big)\cdot\nabla\phi\dx .
\]

\emph{Collecting.} Expanding $\tp^2=(p-p_D)^2$ and $\tn^2=(n-n_D)^2$ and using $(p^2-n^2)q=s\,q^2$, we obtain
\[
\big(\tp^2-\tn^2\big)q=s\,q^2-2\big(pp_D-nn_D\big)q+\big(p_D^2-n_D^2\big)q .
\]
Keeping $\frac{\mu_+\mu_-}{2\epsilon_0}\int_\Om s\,q^2\dx$ on the left and moving the data terms to the right gives \eqref{eq:L2identity}.
\end{proof}

The screening term $\frac{\mu_+\mu_-}{2\epsilon_0}\int_\Om s\,q^2\dx$, which expresses the dissipativity of the Poisson coupling, is nonnegative because $p,n\ge0$. The data terms are bounded as follows.

\begin{lemma}[Data terms]\label{lem:data}
There is a constant $C$ depending only on $\mathcal D$ such that, for a.e.\ $t\in(0,T)$,
\begin{equation}\label{est:data}
    \abs{\Rdata}\le\frac{\theta}{2}\Big(\norm{\nabla\tp}_2^2+\norm{\nabla\tn}_2^2\Big)+C\Big(1+\norm{\nabla\phi}_2^2+\Y\Big).
\end{equation}
\end{lemma}

\begin{proof}
By Young's inequality, since $\nabla p_D\in L^2(\Om)$ and $p_D\in L^\infty(\Om)$, we obtain
\[
\epsilon_+\mu_-\Big|\int_\Om\nabla p_D\cdot\nabla\tp\dx\Big|\le\frac\theta8\norm{\nabla\tp}_2^2+C,\qquad
\mu_+\mu_-\Big|\int_\Om p_D\nabla\tp\cdot\nabla\phi\dx\Big|\le\frac\theta8\norm{\nabla\tp}_2^2+C\norm{\nabla\phi}_2^2,
\]
and
\[
\mu_-\Big|\int_\Om p_D\,\bv\cdot\nabla\tp\dx\Big|\le\mu_-\norm{p_D}_2\norm{\bv}_\infty\norm{\nabla\tp}_2\le\frac\theta8\norm{\nabla\tp}_2^2+\frac{2\mu_-^2}{\theta}\norm{p_D}_2^2\norm{\bv}_\infty^2 ;
\]
the same bounds hold for the $\tn$-terms. Since $\abs{q}\le s$ and $p_D,n_D\in L^\infty(\Om)$,
\begin{gather*}
\frac{\mu_+\mu_-}{\epsilon_0}\Big|\int_\Om(pp_D-nn_D)q\dx\Big|\le C\int_\Om s^2\dx\le C(1+\Y),\\
\frac{\mu_+\mu_-}{2\epsilon_0}\Big|\int_\Om(p_D^2-n_D^2)q\dx\Big|\le C\norm{q}_2\le C(1+\Y).
\end{gather*}
Adding these eight bounds, in which $\norm{\nabla\tp}_2^2$ and $\norm{\nabla\tn}_2^2$ each occur three times with the coefficient $\theta/8$, we obtain
\[
\abs{\Rdata}\le\frac{3\theta}{8}\Big(\norm{\nabla\tp}_2^2+\norm{\nabla\tn}_2^2\Big)+C\Big(1+\norm{\nabla\phi}_2^2+\Y\Big),
\]
which implies the asserted estimate since $3\theta/8\le\theta/2$.
\end{proof}

The source term admits two estimates. The first uses the H\"older exponents of \cite[(4.13)]{GHL26}; the second keeps the field in $L^2$.

\begin{lemma}[The source term]\label{lem:source}
For a.e.\ $t\in(0,T)$ the following two bounds hold.
\begin{enumerate}[leftmargin=2.2em,label=\textup{(\roman*)}]
\item \textup{(Field measured through the densities.)} $\displaystyle \abs{\RF}\le\frac\theta4\Big(\norm{\nabla\tp}_2^2+\norm{\nabla\tn}_2^2\Big)+C\big(1+\Y\big)^2 .$
\item \textup{(Field kept in $L^2$.)} $\displaystyle \abs{\RF}\le\frac\theta4\Big(\norm{\nabla\tp}_2^2+\norm{\nabla\tn}_2^2\Big)+C\Big(1+\norm{\nabla\phi}_2^{2\kappa}\Big)\big(1+\Y\big).$
\end{enumerate}
\end{lemma}

\begin{proof}
By \eqref{eq:Gamma-bound} and $n\ge0$, we get
\begin{equation}\label{eq:RF-start}
\abs{\RF}\le\mu_-\gamma_*(\mu_++\mu_-)\int_\Om n\,|\nabla\phi|\,\big(\abs{\tp}+\abs{\tn}\big)\dx .
\end{equation}
By H\"older's inequality, Lemma~\ref{lem:aux}(i)--(ii) and \eqref{eq:weak-poisson}, we have
\begin{align*}
\int_\Om n|\nabla\phi|(\abs{\tp}+\abs{\tn})\dx&\le\norm{\nabla\phi}_6\norm{n}_3\norm{\abs{\tp}+\abs{\tn}}_2\\
&\le C\big(1+\norm{q}_2\big)\norm{n}_{H^1}\Y^{1/2}\le C\big(1+\Y^{1/2}\big)\norm{n}_{H^1}\Y^{1/2},
\end{align*}
where we use
$$\norm{\nabla\phi}_6\le C\norm{\phi}_{H^2}\le C(\norm{\phi_D}_{H^2}+\norm{\psi}_{H^2})\le C(1+\norm{q}_2)\quad \text{and}\quad \norm{q}_2\le C(1+\Y^{1/2}).$$ 
Young's inequality and $\norm{n}_{H^1}^2\le C(\norm{\nabla\tn}_2^2+1+\Y)$ give (i).

Since $\abs{\tp}+\abs{\tn}\le s+p_D+n_D$, $0\le n\le s$ and $p_D,n_D\in L^\infty(\Om)$, we obtain
\[
\int_\Om n|\nabla\phi|(\abs{\tp}+\abs{\tn})\dx\le\norm{\nabla\phi}_2\norm{ns}_2+C\norm{\nabla\phi}_2\norm{n}_2
\le\norm{\nabla\phi}_2\norm{s}_4^2+C\big(\norm{\nabla\phi}_2^2+1+\Y\big).
\]
If $d=2$, Lemma~\ref{lem:aux}(i) and Young's inequality give
\[
\norm{\nabla\phi}_2\norm{s}_4^2\le C_{GN}\norm{\nabla\phi}_2\norm{s}_2\norm{s}_{H^1}\le\delta\norm{s}_{H^1}^2+C_\delta\norm{\nabla\phi}_2^2\norm{s}_2^2 ,
\]
and if $d=3$, Lemma~\ref{lem:aux}(i) and Young's inequality with exponents $\frac43$ and $4$ give
\[
\norm{\nabla\phi}_2\norm{s}_4^2\le C_{GN}\norm{\nabla\phi}_2\norm{s}_2^{1/2}\norm{s}_{H^1}^{3/2}\le\delta\norm{s}_{H^1}^2+C_\delta\norm{\nabla\phi}_2^4\norm{s}_2^2 .
\]
In both cases 
\[\norm{s}_{H^1}^2\le C(\norm{\nabla\tp}_2^2+\norm{\nabla\tn}_2^2+1+\Y)\quad \text{and}\quad \norm{s}_2^2\le C(1+\Y).\]
Choosing $\delta$ small and using $\norm{\nabla\phi}_2^2\le1+\norm{\nabla\phi}_2^{2\kappa}$ yields (ii).
\end{proof}

Inserting Lemma~\ref{lem:data} and Lemma~\ref{lem:source}(i) into \eqref{eq:L2identity}, and using $\norm{\nabla\phi}_2^2\le C(1+\norm{q}_2^2)\le C(1+\Y)$, we obtain a differential inequality of Riccati type.

\begin{corollary}[Riccati-type inequality]\label{cor:riccati}
There is $C_0=C_0(\mathcal D)$ such that
\[
\frac{d\Y}{dt}+\frac\theta2\Big(\norm{\nabla\tp}_2^2+\norm{\nabla\tn}_2^2\Big)+\frac{\mu_+\mu_-}{\epsilon_0}\int_\Om s\,q^2\dx\le C_0\big(1+\Y\big)^2\qquad\text{a.e.\ in }(0,T).
\]
Consequently $1+\Y(t)\le\big(1+\Y(0)\big)\big(1-C_0(1+\Y(0))\,t\big)^{-1}$ for $0\le t<T$ with $C_0(1+\Y(0))\,t<1$.
\end{corollary}

The bound of Corollary~\ref{cor:riccati} gives no information for $t\ge1/\big(C_0(1+\Y(0))\big)$. Corollary~\ref{cor:riccati} is the differential form, for weak solutions of \eqref{eq:sys}--\eqref{eq:bc}, of the integral inequality \cite[(4.14)]{GHL26}, written as \eqref{eq:GHL-riccati}, which concerns the truncated problem. There, the Bihari inequality \cite[Proposition~2.3]{GHL26} with the nonlinearity $z\mapsto z+z^2$ gives the bound \cite[(4.15)]{GHL26}, which holds up to the blow-up time of the comparison equation, and the existence time $T_1$ of \cite[Theorem~1.1]{GHL26} is fixed accordingly (Section~\ref{ssec:model}). In both cases it is the quadratic term that confines the estimate to a finite time interval.

\begin{remark}[Where the quadratic term comes from]\label{rem:origin}
In \cite[(4.14)]{GHL26} the quadratic term comes from the ionization contribution of \cite[(3.15)]{GHL26}, written here without the superscript $M$,
\[
J_1=\mu_-\alpha_1\int_0^t\!\!\int_\Om e^{-\alpha_2/|\nabla\phi|}|\nabla\phi|\,n\,\big(\mu_-\tp+\mu_+\tn\big)\dx\dr ,
\]
and in Corollary~\ref{cor:riccati} it comes from $\RF$, through Lemma~\ref{lem:source}(i). In both cases the spatial integrand is estimated, in \cite[(4.13)]{GHL26} and in the proof of Lemma~\ref{lem:source}(i), by H\"older's inequality with exponents $6,3,2$, schematically
\begin{equation}\label{eq:J1-holder}
\int_\Om|\nabla\phi|\,n\,(p+n)\dx\le\norm{\nabla\phi}_6\,\norm{n}_3\,\norm{p+n}_2 .
\end{equation}
The first factor is bounded through the Poisson equation and Lemma~\ref{lem:aux}(i)--(ii), $\norm{\nabla\phi}_6\le C(1+\norm{q}_2)$, and the second through Lemma~\ref{lem:aux}(i), $\norm{n}_3\le C\norm{n}_{H^1}$. Young's inequality then bounds the right-hand side of \eqref{eq:J1-holder} by
\[
\delta\norm{n}_{H^1}^2+C_\delta\big(1+\norm{q}_2\big)^2\norm{p+n}_2^2 ,
\]
and $\delta\norm{n}_{H^1}^2$ is absorbed by the dissipation, in \eqref{eq:L2identity} or on the left of \cite[(4.11)]{GHL26}, up to $C(1+\Y)$. Since $\norm{q}_2,\norm{p+n}_2\le\norm{p}_2+\norm{n}_2$, each of the two remaining factors is bounded by $C(1+\Y)$, and their product by $C(1+\Y)^2$.

Lemma~\ref{lem:source}(ii) estimates the same integrand with exponents $2,4,4$: it keeps $\norm{\nabla\phi}_2$ instead of bounding $\norm{\nabla\phi}_6$ by $C(1+\norm{q}_2)$, and it bounds $\norm{p+n}_4^2$ by interpolation between $\norm{p+n}_2$ and $\norm{p+n}_{H^1}$. The field then enters as a coefficient rather than through a power of $\Y$.
\end{remark}

Using Lemma~\ref{lem:source}(ii) in place of Lemma~\ref{lem:source}(i) we obtain an inequality that is linear in $\Y$, with the field in the coefficient.

\begin{proposition}[Field-controlled energy inequality]\label{prop:field-controlled}
There is a constant $C_1$ depending only on $\mathcal D$ such that
\begin{equation}\label{eq:field-controlled}
\frac{d\Y}{dt}+\frac\theta2\Big(\norm{\nabla\tp}_2^2+\norm{\nabla\tn}_2^2\Big)+\frac{\mu_+\mu_-}{\epsilon_0}\int_\Om s\,q^2\dx\le a(t)\,\big(1+\Y\big)
\end{equation}
holds a.e. in $(0, T)$,  where
\begin{equation}\label{eq:coefficient}
a(t):=C_1\Big(1+\norm{\nabla\phi(t)}_2^{2\kappa}\Big),
\end{equation} 
and therefore, for $0\le t<T$,
\begin{equation}\label{eq:field-gronwall}
1+\Y(t)+\frac\theta2\int_0^t\Big(\norm{\nabla\tp}_2^2+\norm{\nabla\tn}_2^2\Big)\dr\le\big(1+\Y(0)\big)\exp\Big(\int_0^ta(r)\dr\Big).
\end{equation}
\end{proposition}

\begin{proof}
Inserting Lemma~\ref{lem:data} and Lemma~\ref{lem:source}(ii) into \eqref{eq:L2identity}, multiplying by $2$, and using the inequality $\norm{\nabla\phi}_2^2\le1+\norm{\nabla\phi}_2^{2\kappa}$ for $\kappa\ge1$ gives, for a.e.\ $t\in(0,T)$,
\[\begin{split}
&\frac{d\Y}{dt}+2\epsilon_+\mu_-\norm{\nabla\tp}_2^2+2\epsilon_-\mu_+\norm{\nabla\tn}_2^2+\frac{\mu_+\mu_-}{\epsilon_0}\int_\Om s\,q^2\dx\\
\le&\frac{3\theta}{2}\Big(\norm{\nabla\tp}_2^2+\norm{\nabla\tn}_2^2\Big)+C\Big(1+\norm{\nabla\phi}_2^{2\kappa}\Big)\big(1+\Y\big).
\end{split}\]
Since $\epsilon_+\mu_-,\epsilon_-\mu_+\ge\theta$ by \eqref{eq:Y}, this implies \eqref{eq:field-controlled} with a constant $C_1$ depending only on $\mathcal D$.

Put $A(t):=\int_0^ta(r)\dr$, which is Lipschitz continuous and nondecreasing on $[0,T]$ since $a\ge0$ and $a\in L^\infty(0,T)$ by Definition~\ref{def:weak}(i). As $\Y$ is absolutely continuous on $[0,T]$ (proof of Lemma~\ref{lem:L2identity}), so is $e^{-A}(1+\Y)$. Multiplying \eqref{eq:field-controlled} by $e^{-A(t)}$ and dropping the nonnegative screening term $\frac{\mu_+\mu_-}{\epsilon_0}\int_\Om s\,q^2\dx$, we obtain, for a.e.\ $t\in(0,T)$,
\[
\frac{d}{dt}\Big(e^{-A(t)}\big(1+\Y(t)\big)\Big)+\frac\theta2\,e^{-A(t)}\Big(\norm{\nabla\tp(t)}_2^2+\norm{\nabla\tn(t)}_2^2\Big)\le0 .
\]
Integrating this inequality over $(0,t)$ and using $e^{-A(r)}\ge e^{-A(t)}$ for $0\le r\le t$ yields
\[
e^{-A(t)}\bigg(1+\Y(t)+\frac\theta2\int_0^t\Big(\norm{\nabla\tp}_2^2+\norm{\nabla\tn}_2^2\Big)\dr\bigg)\le1+\Y(0),
\]
which is equivalent to \eqref{eq:field-gronwall}.
\end{proof}

\begin{proof}[Proof of Theorem~\ref{thm:A}]
Let $T_{\max}<\infty$ and suppose, for contradiction, that
\[
I:=\int_0^{T_{\max}}\norm{\nabla\phi}_2^{2\kappa}\dt<\infty .
\]
Given $t\in(0,T_{\max})$, pick $T'\in(t,T_{\max})$. The maximal solution is a weak solution on $(0,T')$, so \eqref{eq:field-gronwall} holds at $t$ with the constant $C_1$ of \eqref{eq:coefficient}, which depends only on $\mathcal D$. Since $$\int_0^ta(r)\dr\le C_1(T_{\max}+I),$$ we obtain 
\begin{equation}\label{eq:uniform-K}
1+\Y(t)\le\big(1+\Y(0)\big)\exp\big(C_1(T_{\max}+I)\big)=:K\qquad\text{for all }t\in(0,T_{\max}).
\end{equation}
Hence, by the bound $\norm{p}_2^2+\norm{n}_2^2\le C(1+\Y)$ stated after \eqref{eq:Y},
\[
\sup_{0\le t<T_{\max}}\big(\norm{p(t)}_2+\norm{n(t)}_2\big)\le CK^{1/2}<\infty,
\]
which contradicts Lemma~\ref{lem:alternative}. This proves the first assertion \eqref{thm:A_1st_conclusion}.

For the second, write $\phi=\phi_D+\psi$. By \eqref{eq:weak-poisson}, the weak form of \eqref{eq:poisson}, $\psi(t)$ is the solution of the Poisson problem of Lemma~\ref{lem:aux}(ii) with $f=q(t)$, so $\norm{\nabla\psi}_2\le\epsilon_0^{-1}\norm{q}_{\HDd}$ and
\begin{equation}\label{est:phi_2kappa}
\norm{\nabla\phi}_2^{2\kappa}\le\left(\norm{\nabla\phi_D}_2+\epsilon_0^{-1}\norm{q}_{\HDd}\right)^{2\kappa}\le C\left(1+\norm{q}_{\HDd}^{2\kappa}\right).
\end{equation}
Integrating \eqref{est:phi_2kappa} over $(0,T_{\max})$ and using \eqref{thm:A_1st_conclusion} gives
\[
\infty=\int_0^{T_{\max}}\norm{\nabla\phi}_2^{2\kappa}\dt\le C\left(T_{\max}+\int_0^{T_{\max}}\norm{q}_{\HDd}^{2\kappa}\dt\right),
\]
and since $T_{\max}<\infty$, the integral $\int_0^{T_{\max}}\norm{q}_{\HDd}^{2\kappa}\dt$ must be infinite, which is \eqref{thm:A_2nd_conclusion}.

For the third \eqref{thm:A_3rd_conclusion}, suppose that $$L:=\limsup_{t\to T_{\max}}\norm{q(t)}_{\HDd}<\infty.$$ 
Then there is $0<S<T_{\max}$ with $\norm{q(t)}_{\HDd}\le L+1$ for all $t\in(S,T_{\max})$. Moreover, $q\in C([0,S],L^2(\Om))$ (Section~\ref{ssec:notation}) and $\norm{f}_{\HDd}\le C_P\norm{f}_2$, so $\norm{q(t)}_{\HDd}\le C_P\max_{[0,S]}\norm{q}_2<\infty$ for all $t\in[0,S]$ as well. Hence $M_0:=\sup_{0<t<T_{\max}}\norm{q(t)}_{\HDd}<\infty$, and so $$\int_0^{T_{\max}}\norm{q}_{\HDd}^{2\kappa}\dt\le T_{\max}M_0^{2\kappa}<\infty,$$ which contradicts \eqref{thm:A_2nd_conclusion}. This proves \eqref{thm:A_3rd_conclusion}.
\end{proof}

\begin{remark}[Blow-up must be a space-charge singularity]
By \eqref{eq:weak-poisson} and Lemma~\ref{lem:aux}(ii), the electrostatic energy $\ee$ of \eqref{eq:E-J} satisfies
\begin{equation}\label{energy=Poisson}
\ee=\frac{\epsilon_0}{2}\norm{\nabla\psi}_2^2=\frac{1}{2\epsilon_0}\norm{p-n}_{\HDd}^2 .
\end{equation}
Theorem~\ref{thm:A} therefore excludes blow-up while $\ee$ stays bounded. Equivalently, no ``quasi-neutral'' blow-up, with $T_{\max}<\infty$ and $p-n$ bounded in $\HDd$ on $(0,T_{\max})$, is possible. Since $$\norm{p-n}_{\HDd}\le C_P\big(\norm{p}_2+\norm{n}_2\big),$$ the conclusion of Theorem~\ref{thm:A} implies that of Lemma~\ref{lem:alternative}, but not conversely. The proof of Theorem~\ref{thm:A} rests on two facts: by Lemma~\ref{lem:L2identity}, drift and convection contribute, apart from data terms, only the nonnegative screening term; and the source $\Gamma(|\nabla\phi|)\,n$ is linear in $n$, with a rate bounded by the field through \eqref{eq:Gamma-bound}.
\end{remark}

\section{Global existence: proof of Theorem~\ref{thm:B}}\label{sec:global}

By Theorem~\ref{thm:A} and the inequality $\norm{\nabla\phi}_2\le\norm{\nabla\phi_D}_2+(2\ee/\epsilon_0)^{1/2}$, global existence follows from a bound on the electrostatic energy $\ee$ of the space charge on bounded time intervals. We bound $\ee$ together with $m=\norm{p-p_D}_1+\norm{n-n_D}_1$, the $L^1$-distance of the densities from the data $p_D,n_D$, which controls the number of carriers $\int_\Om s\dx$ (see \eqref{eq:mass}). Section~\ref{sec:electro} derives a balance law for $\ee$, Section~\ref{sec:L1} an $L^1$ estimate for $m$, and Section~\ref{sec:proofB} closes the two by a linear Gronwall inequality. Throughout Sections~\ref{sec:electro} and~\ref{sec:L1}, $(\phi,p,n)$ is a weak solution on $(0,T)$ under Assumption~\ref{ass:main}, with the abbreviations \eqref{eq:abbrev}.

\subsection{The electrostatic energy balance}\label{sec:electro}

Let $\psi=\phi-\phi_D$ as in Definition~\ref{def:weak}(ii), and define the electrostatic energy of the space charge and the Joule heating by
\begin{equation}\label{eq:E-J}
\ee(t):=\frac{\epsilon_0}{2}\int_\Om\abs{\nabla\psi(t)}^2\dx,\qquad
\Jou(t):=\int_\Om g(t)\abs{\nabla\psi(t)}^2\dx=\int_\Om\big(\mu_+p+\mu_-n\big)\abs{\nabla\psi}^2\dx .
\end{equation}
Since $p,n\ge0$, we have $g\ge0$ and $\Jou\ge0$. By \eqref{eq:psi-regularity}, $\ee$ is absolutely continuous on $[0,T]$. We write $q_D:=p_D-n_D$.

\begin{lemma}[Electrostatic energy balance]\label{lem:energy-identity}
For a.e.\ $t\in(0,T)$,
\begin{equation}\label{eq:energy-identity}
\begin{aligned}
\frac{d\ee}{dt}+\Jou+\frac{\epsilon_++\epsilon_-}{2\epsilon_0}\norm{q}_2^2
&=-\frac{\epsilon_+-\epsilon_-}{2\epsilon_0}\int_\Om\big(p^2-n^2\big)\dx
-\int_\Om g\nabla\phi_D\cdot\nabla\psi\dx\\
&\quad-\int_\Om\nabla\big(\epsilon_+p_D-\epsilon_-n_D\big)\cdot\nabla\psi\dx
+\frac1{\epsilon_0}\int_\Om\big(\epsilon_+p_D-\epsilon_-n_D\big)q\dx\\
&\quad+\int_\Om q\bv\cdot\nabla\psi\dx .
\end{aligned}
\end{equation}
\end{lemma}

\begin{proof}
By \eqref{eq:psi-regularity} with $\chi=\psi(t)$,
\[
\frac{d\ee}{dt}=\epsilon_0\int_\Om\nabla\psi\cdot\nabla\partial_t\psi\dx=\pair{\partial_tq,\psi}=\pair{\partial_tp,\psi}-\pair{\partial_tn,\psi}.
\]
We now take $f=\psi(t)\in\HD$ in \eqref{eq:weak-p} and in \eqref{eq:weak-n} and subtract. The source terms $\int_\Om F\psi\dx$ cancel identically, and we obtain
\begin{equation}\label{eq:dE-raw}
\frac{d\ee}{dt}=-\int_\Om\big(\epsilon_+\nabla p-\epsilon_-\nabla n\big)\cdot\nabla\psi\dx-\int_\Om g\nabla\phi\cdot\nabla\psi\dx+\int_\Om q\bv\cdot\nabla\psi\dx .
\end{equation}
Since $\nabla\phi=\nabla\psi+\nabla\phi_D$, the second term on the right-hand side of \eqref{eq:dE-raw} becomes  
\begin{equation}\label{eq:dE-raw_2nd_term}-\int_\Om g\nabla\phi\cdot\nabla\psi\dx=-\Jou-\int_\Om g\nabla\phi_D\cdot\nabla\psi\dx.\end{equation} For the first term on the right-hand side of \eqref{eq:dE-raw} we write $$\epsilon_+\nabla p-\epsilon_-\nabla n=\nabla(\epsilon_+\tp-\epsilon_-\tn)+\nabla(\epsilon_+p_D-\epsilon_-n_D).$$
Since $\epsilon_+\tp-\epsilon_-\tn\in\HD$, \eqref{eq:weak-poisson} with $\chi=\epsilon_+\tp-\epsilon_-\tn$ gives
\begin{align*}
-\int_\Om\nabla\big(\epsilon_+\tp-\epsilon_-\tn\big)\cdot\nabla\psi\dx&=-\frac1{\epsilon_0}\int_\Om q\big(\epsilon_+\tp-\epsilon_-\tn\big)\dx\\
&=-\frac1{\epsilon_0}\int_\Om q\big(\epsilon_+p-\epsilon_-n\big)\dx+\frac1{\epsilon_0}\int_\Om q\big(\epsilon_+p_D-\epsilon_-n_D\big)\dx .
\end{align*}
Finally, $\epsilon_+p-\epsilon_-n=\frac{\epsilon_++\epsilon_-}{2}\,q+\frac{\epsilon_+-\epsilon_-}{2}\,s$ and $sq=p^2-n^2$, so that the first term on the right-hand side of \eqref{eq:dE-raw} becomes
\begin{multline*}
-\int_\Om\big(\epsilon_+\nabla p-\epsilon_-\nabla n\big)\cdot\nabla\psi\dx=-\frac{\epsilon_++\epsilon_-}{2\epsilon_0}\norm{q}_2^2-\frac{\epsilon_+-\epsilon_-}{2\epsilon_0}\int_\Om\big(p^2-n^2\big)\dx\\
-\int_\Om\nabla\big(\epsilon_+p_D-\epsilon_-n_D\big)\cdot\nabla\psi\dx+\frac1{\epsilon_0}\int_\Om\big(\epsilon_+p_D-\epsilon_-n_D\big)\,q\dx .
\end{multline*}
Inserting this and the expression \eqref{eq:dE-raw_2nd_term} for the second term into \eqref{eq:dE-raw}, and moving $\Jou$ and $\frac{\epsilon_++\epsilon_-}{2\epsilon_0}\norm{q}_2^2$ to the left-hand side, gives \eqref{eq:energy-identity}.
\end{proof}

\begin{remark}[Interpretation]\label{rem:interpretation}
Identity \eqref{eq:energy-identity} is the energy balance of the electrostatic field generated by the space charge. By \eqref{eq:dE-raw},
\[
\frac{d\ee}{dt}+\int_\Om\bj\cdot(-\nabla\psi)\dx=0,\qquad \bj:=q\bv-\epsilon_+\nabla p+\epsilon_-\nabla n-g\nabla\phi ,
\]
where $\bj$ is the charge current, which satisfies $\partial_tq+\nabla\cdot\bj=0$ by \eqref{eq:charge-eq} and $\nabla\cdot\bv=0$, and the integral is the rate of work, that is, the work per unit time, done by the space-charge field $-\nabla\psi$ on $\bj$. Hence $\ee$ decreases while the rate of work is positive and increases while it is negative; in general, the rate of work has no sign. Since $\nabla\phi=\nabla\psi+\nabla\phi_D$, the drift current $-g\nabla\phi$ splits into a part $-g\nabla\psi$ driven by the space-charge field itself and a part $-g\nabla\phi_D$ driven by the applied field. The rate of work on the first part is $\int_\Om(-g\nabla\psi)\cdot(-\nabla\psi)\dx=\Jou$, the Joule heating, which is nonnegative since $g\ge0$. The rate of work on the second part is $\int_\Om g\,\nabla\phi_D\cdot\nabla\psi\dx$; we call its negative, $-\int_\Om g\,\nabla\phi_D\cdot\nabla\psi\dx$, the battery term: through it the voltage source acts on $\ee$. By the proof of Lemma~\ref{lem:energy-identity}, the rate of work on the diffusion current $-\epsilon_+\nabla p+\epsilon_-\nabla n$ is
\begin{multline*}
\frac{\epsilon_++\epsilon_-}{2\epsilon_0}\norm{q}_2^2+\frac{\epsilon_+-\epsilon_-}{2\epsilon_0}\int_\Om\big(p^2-n^2\big)\dx\\
+\int_\Om\nabla\big(\epsilon_+p_D-\epsilon_-n_D\big)\cdot\nabla\psi\dx-\frac1{\epsilon_0}\int_\Om\big(\epsilon_+p_D-\epsilon_-n_D\big)\,q\dx .
\end{multline*}
The first term, which we call the dissipative term, is nonnegative. We call the negative of the second term, $-\frac{\epsilon_+-\epsilon_-}{2\epsilon_0}\int_\Om(p^2-n^2)\dx$, the ambipolar term; it describes the exchange of energy between the field and the carriers by differential (ambipolar) diffusion, has no sign and vanishes under (H1). The last two terms are data terms: they are linear in $\psi$ and $q$, with coefficients given by $p_D$ and $n_D$. Finally, the rate of work on the convection current $q\bv$ is $-\int_\Om q\,\bv\cdot\nabla\psi\dx$. The source $F$ does not appear at all: it creates or removes carriers in pairs and so does not change the charge.
\end{remark}

\begin{corollary}[Energy inequality under \textup{(H1)}--\textup{(H2)}]\label{cor:energy-ineq}
Assume in addition \textup{(H1)} and \textup{(H2)}. Then there is $C_2=C_2(\mathcal D)$ such that for a.e.\ $t\in(0,T)$
\begin{equation}\label{eq:energy-ineq}
\frac{d\ee}{dt}+\frac12\Jou+\frac{\epsilon}{2\epsilon_0}\norm{q}_2^2\le C_2\Big(1+\ee+\int_\Om s\dx\Big).
\end{equation}
\end{corollary}

\begin{proof}
Under (H1) the first term on the right of \eqref{eq:energy-identity} vanishes and the coefficient of $\norm{q}_2^2$ on the left is $\epsilon/\epsilon_0$. We bound the remaining four terms on the right. For the battery term, (H2), the Cauchy--Schwarz inequality with respect to the measure $g\dx$, Young's inequality and $g\le\max\{\mu_+,\mu_-\}s$ give
\begin{equation}\label{eq:battery-bound}
\begin{split}\Big|\int_\Om g\,\nabla\phi_D\cdot\nabla\psi\dx\Big|&\le\norm{\nabla\phi_D}_\infty\Big(\int_\Om g\dx\Big)^{1/2}\Jou^{1/2}\\
&\le\frac12\Jou+\frac12\norm{\nabla\phi_D}_\infty^2\max\{\mu_+,\mu_-\}\int_\Om s\dx .
\end{split}\end{equation}
Under (H1) the two data terms on the right of \eqref{eq:energy-identity} are 
\[-\int_\Om\nabla\big(\epsilon_+p_D-\epsilon_-n_D\big)\cdot\nabla\psi\dx=-\epsilon\int_\Om\nabla q_D\cdot\nabla\psi\dx, 
\quad   \frac1{\epsilon_0}\int_\Om\big(\epsilon_+p_D-\epsilon_-n_D\big)q\dx=\frac{\epsilon}{\epsilon_0}\int_\Om q_Dq\dx.\] 
Since $\norm{\nabla\psi}_2=(2\ee/\epsilon_0)^{1/2}\le C(1+\ee)$, the Cauchy--Schwarz and Young inequalities give
\begin{equation}\label{est:2-3terms_energy_eqy}\begin{split}
\epsilon\Big|\int_\Om\nabla q_D\cdot\nabla\psi\dx\Big|\le\epsilon\norm{\nabla q_D}_2\norm{\nabla\psi}_2\le C(1+\ee),\\
\frac{\epsilon}{\epsilon_0}\Big|\int_\Om q_Dq\dx\Big|\le\frac{\epsilon}{4\epsilon_0}\norm{q}_2^2+\frac{\epsilon}{\epsilon_0}\norm{q_D}_2^2,\qquad\qquad 
\end{split}\end{equation}
and for the convection term 
\begin{equation}\label{est:convection}
\Big|\int_\Om q\bv\cdot\nabla\psi\dx\Big|\le\norm{\bv}_\infty\norm{q}_2\norm{\nabla\psi}_2\le\frac{\epsilon}{4\epsilon_0}\norm{q}_2^2+\frac{\epsilon_0}{\epsilon}\norm{\bv}_\infty^2\norm{\nabla\psi}_2^2\le\frac{\epsilon}{4\epsilon_0}\norm{q}_2^2+C\ee .
\end{equation}
Inserting \eqref{eq:battery-bound}, \eqref{est:2-3terms_energy_eqy} and \eqref{est:convection} into \eqref{eq:energy-identity} and using (H1) gives
\[\begin{split}\frac{d\ee}{dt}+\Jou+\frac{\epsilon}{\epsilon_0}\norm{q}_2^2\le \frac12\Jou+\frac12\norm{\nabla\phi_D}_\infty^2\max\{\mu_+,\mu_-\}\int_\Om s\dx\\
+ C(1+\ee)+\frac{\epsilon}{2\epsilon_0}\norm{q}_2^2+\frac{\epsilon}{\epsilon_0}\norm{q_D}_2^2+C\ee.\quad\end{split}\]
This concludes the proof of \eqref{eq:energy-ineq}.
\end{proof}

Apart from $\ee$ and the data, the right-hand side of \eqref{eq:energy-ineq} involves the densities only through the number of carriers $\int_\Om s\dx$, which Section~\ref{sec:L1} estimates, and not through their $L^2$-norms. In the $L^2$ estimate \cite[(4.14)]{GHL26}, written as \eqref{eq:GHL-riccati}, the densities enter the right-hand side through $Y$ and $Y^2$.

\begin{remark}[Why the free energy does not close the argument]\label{rem:freeenergy}
For $F=0$ and $\bv=0$, \eqref{eq:sys} is the bipolar drift--diffusion (van Roosbroeck) system, for which global well-posedness is classical under the assumptions of \cite{Gaj85,GG86} (Section~\ref{ssec:background}). With the weights $\epsilon_\pm/\mu_\pm$, its entropy and free energy are
\[
\mathcal S:=\frac{\epsilon_+}{\mu_+}\int_\Om\big(p\log p-p\big)\dx+\frac{\epsilon_-}{\mu_-}\int_\Om\big(n\log n-n\big)\dx,\qquad
\mathcal F:=\mathcal S+\ee .
\]
Formally, for $\bv=0$, $\phi_D=0$ (so that $\psi=\phi$) and no-flux conditions for $p$ and $n$, testing \eqref{eq:p} with $\frac{\epsilon_+}{\mu_+}\log p$ and \eqref{eq:n} with $\frac{\epsilon_-}{\mu_-}\log n$ gives, with \eqref{eq:dE-raw},
\begin{align*}
\frac{d\mathcal S}{dt}&=-\frac{\epsilon_+^2}{\mu_+}\int_\Om\frac{|\nabla p|^2}{p}\dx-\frac{\epsilon_-^2}{\mu_-}\int_\Om\frac{|\nabla n|^2}{n}\dx-\mathcal C+\int_\Om F\Big(\tfrac{\epsilon_+}{\mu_+}\log p+\tfrac{\epsilon_-}{\mu_-}\log n\Big)\dx,\\
\frac{d\ee}{dt}&=-\mathcal C-\Jou,\qquad
\mathcal C:=\int_\Om\big(\epsilon_+\nabla p-\epsilon_-\nabla n\big)\cdot\nabla\phi\dx .
\end{align*}
In the proof of Lemma~\ref{lem:energy-identity}, the cross term $-\mathcal C$ (the first term on the right of \eqref{eq:dE-raw} with $\psi=\phi$) is rewritten through the Poisson equation, and the ambipolar term $-\frac{\epsilon_+-\epsilon_-}{2\epsilon_0}\int_\Om(p^2-n^2)\dx$ that results has no sign; it vanishes under (H1). In the computation of $d\mathcal F/dt$ the cross term need not be rewritten: $-\mathcal C$ appears in both $d\mathcal S/dt$ and $d\ee/dt$, and with
\begin{align*}
\mathcal I&:=\mu_+\int_\Om p\,\Big|\nabla\Big(\tfrac{\epsilon_+}{\mu_+}\log p+\phi\Big)\Big|^2\dx+\mu_-\int_\Om n\,\Big|\nabla\Big(\tfrac{\epsilon_-}{\mu_-}\log n-\phi\Big)\Big|^2\dx\\
&=\mathcal{K}+2\,\mathcal C+\Jou,\quad \mathcal{K}:=\frac{\epsilon_+^2}{\mu_+}\int_\Om\frac{|\nabla p|^2}{p}\dx+\frac{\epsilon_-^2}{\mu_-}\int_\Om\frac{|\nabla n|^2}{n}\dx,
\end{align*}
it follows that, for arbitrary $\epsilon_\pm,\mu_\pm$,
\[
\frac{d\mathcal F}{dt}+\mathcal I=\int_\Om F\Big(\tfrac{\epsilon_+}{\mu_+}\log p+\tfrac{\epsilon_-}{\mu_-}\log n\Big)\dx,\qquad\mathcal I\ge0 .
\]
No ambipolar term arises, and (H1) is not needed. For $\phi_D\neq0$, the same computation, with $\psi$ in place of $\phi$ in $\mathcal C$ and $\mathcal I$, adds the term $-\int_\Om\big(g\nabla\psi+\epsilon_+\nabla p-\epsilon_-\nabla n\big)\cdot\nabla\phi_D\dx$ to the right-hand side of the identity for $d\mathcal F/dt$; it contains the battery term of Remark~\ref{rem:interpretation}. However, $\mathcal I$ does not control $\Jou$, that is, no inequality $\Jou\le C\,\mathcal I$ can hold: $\mathcal I=0$ on every state with $\nabla\phi=\frac{\epsilon_-}{\mu_-}\nabla\log n=-\frac{\epsilon_+}{\mu_+}\nabla\log p$, for instance on the Poisson--Boltzmann equilibria $p=c_+e^{-\mu_+\phi/\epsilon_+}$, $n=c_-e^{\mu_-\phi/\epsilon_-}$ with $c_+\neq c_-$, whereas on such a state $\Jou>0$ unless $p=n$ (if $\mathcal J=0$, then $\nabla\phi=0$ a.e., since $g=\mu_+p+\mu_-n>0$ on such a state; so $\phi$
is constant, the constant is $0$ because $\phi=0$ on $\GD$, and $q=-\epsilon_0\Delta\phi=0$); hence $\mathcal I=0<\Jou$ for $c_+\neq c_-$. But $\Jou$ bounds the battery term, by \eqref{eq:battery-bound} in the proof of Corollary~\ref{cor:energy-ineq}, and the source of the $L^1$ estimate of Lemma~\ref{lem:L1}, by \eqref{eq:source-Joule}. Since $\mathcal I$ does not control $\Jou$, the factors $\Jou^{1/2}$ on the right-hand sides of \eqref{eq:battery-bound} and \eqref{eq:source-Joule} cannot be absorbed by $\mathcal I$ (the only nonnegative term on the left of the identity for $d\mathcal F/dt$), nor by $\mathcal I$ together with a multiple of $1+\mathcal F+\int_\Om s\dx$, as a linear Gronwall argument would allow: for $c_\pm=\varrho\,\bar c_\pm$ with fixed $\bar c_+\neq\bar c_-$ and $\varrho\to\infty$, the field of the Poisson--Boltzmann state is confined to layers of width of order $\varrho^{-1/2}$ at $\GD$, so that $\Jou$ is of order $\varrho^{3/2}$, whereas $\mathcal F$ is of order $\varrho\log\varrho$ and $\int_\Om s\dx$ of order $\varrho$.

Moreover, the source term $\int_\Om F\big(\tfrac{\epsilon_+}{\mu_+}\log p+\tfrac{\epsilon_-}{\mu_-}\log n\big)\dx$ has no sign, and the natural estimate of its electron part, by the Cauchy--Schwarz inequality with respect to the measure $n\dx$, calls for $\int_\Om n(\log n)^2\dx$, which $\mathcal F$ does not control. Writing $\ee+\lambda\mathcal S$ with $\lambda\neq1$ does not help: its dissipation $\lambda\mathcal K+(1+\lambda)\mathcal C+\Jou$ is a sum of squares only for $\lambda=1$; otherwise a nonzero multiple of $\mathcal C$ is left over, which produces, as in Lemma~\ref{lem:energy-identity}, a nonzero multiple of the ambipolar term. The electrostatic energy alone keeps $\Jou$ explicit; in return it carries the ambipolar
term of \eqref{eq:energy-identity}, and hence requires (H1) in Corollary~\ref{cor:energy-ineq}.
\end{remark}

\subsection{An \texorpdfstring{$L^1$}{L1} estimate for the carrier densities}\label{sec:L1}

The weak formulation gives no equation for the number of carriers $\int_\Om s\dx$, because the test function $f\equiv1$ does not belong to $\HD$ and the carrier flux through the electrodes is not controlled. We therefore estimate the $L^1$-norms of the differences $\tp=p-p_D$ and $\tn=n-n_D$, which vanish on $\GD$, by testing with a regularized sign of $\tp$ and $\tn$. The Dirichlet data are thereby absorbed, and the boundary integral over $\GN$ produced by the diffusion term vanishes by the Neumann compatibility in (A3). Let
\begin{equation}\label{eq:mass}
m(t):=\int_\Om\big(\abs{\tp(t)}+\abs{\tn(t)}\big)\dx,\qquad\text{so that}\qquad
\int_\Om s(t)\dx\le m(t)+\norm{p_D}_1+\norm{n_D}_1 .
\end{equation}
Since $\tp,\tn\in C([0,T],L^2(\Om))$ and $\Om$ is bounded, $m$ is continuous on $[0,T]$. We call $m$ the mass.

\begin{lemma}[$L^1$ estimate]\label{lem:L1}
There is a constant $C_3$ depending only on $\mathcal D$ such that for all $0\le t_1\le t_2\le T$
\begin{equation}\label{eq:L1}
m(t_2)\le m(t_1)+\int_{t_1}^{t_2}\Big[C_3\Big(1+\norm{\nabla\phi}_2+\norm{q}_2\Big)+2\mu_-\gamma_*\int_\Om n\,\abs{\nabla\phi}\dx\Big]\dt .
\end{equation}
\end{lemma}

\begin{proof}
For $\zeta>0$ let
\begin{gather*}
\Sigma_\zeta(u):=\sqrt{u^2+\zeta^2}-\zeta,\qquad \Sigma_\zeta'(u)=\frac{u}{\sqrt{u^2+\zeta^2}},\qquad \Sigma_\zeta''(u)=\frac{\zeta^2}{(u^2+\zeta^2)^{3/2}},\\
S_\zeta(u):=\int_0^u\xi\,\Sigma_\zeta''(\xi)\,d\xi=\zeta-\frac{\zeta^2}{\sqrt{u^2+\zeta^2}} .
\end{gather*}
Then $\Sigma_\zeta$ is convex, $\Sigma_\zeta(0)=\Sigma_\zeta'(0)=0$, $\abs{\Sigma_\zeta'}\le1$, $\Sigma_\zeta'$ is Lipschitz, $0\le S_\zeta<\zeta$, $S_\zeta$ is Lipschitz with $S_\zeta(0)=0$, and $\Sigma_\zeta(u)\rightarrow\abs{u}$ as $\zeta\rightarrow0$.

\smallskip\noindent\emph{The ion equation.} By Definition~\ref{def:weak}(i), $p\in H^1((0,T),\HDd)$ and $\tp\in L^2((0,T),\HD)$; since $p_D\in L^2(\Om)\subset\HDd$ does not depend on $t$, also $\tp\in H^1((0,T),\HDd)$ with $\partial_t\tp=\partial_tp$. Hence Lemma~\ref{lem:aux}(iii) applies to $u=\tp$ and $\Sigma=\Sigma_\zeta$: the function $t\mapsto\int_\Om\Sigma_\zeta(\tp(t))\dx$ is absolutely continuous on $[0,T]$, with derivative $\pair{\partial_tp,\Sigma_\zeta'(\tp)}$ for a.e.\ $t$. Taking $f=\Sigma_\zeta'(\tp(t))\in\HD\cap L^\infty(\Om)$ in \eqref{eq:weak-p} we obtain, for a.e.\ $t$,
\begin{equation}\label{eq:L1-start}
\frac{d}{dt}\int_\Om\Sigma_\zeta(\tp)\dx=-\int_\Om\big(\epsilon_+\nabla p+\mu_+p\nabla\phi-p\bv\big)\cdot\nabla\Sigma_\zeta'(\tp)\dx+\int_\Om F\Sigma_\zeta'(\tp)\dx.
\end{equation}
The right-hand side of \eqref{eq:L1-start} is written as $I_1+I_2+I_3+I_4$, where 
\begin{align*}
    I_1&:=-\int_\Om\epsilon_+\nabla p\cdot\nabla\Sigma_\zeta'(\tp)\dx, & I_2&:=-\int_\Om\mu_+p\nabla\phi\cdot\nabla\Sigma_\zeta'(\tp)\dx,\\
    I_3&:=\int_\Om p\bv\cdot\nabla\Sigma_\zeta'(\tp)\dx, & I_4&:=\int_\Om F\Sigma_\zeta'(\tp)\dx.
\end{align*}
By the chain rule, $\nabla\Sigma_\zeta'(\tp)=\Sigma_\zeta''(\tp)\nabla\tp$ and $\tp\nabla\Sigma_\zeta'(\tp)=\nabla S_\zeta(\tp)$, with $S_\zeta(\tp)\in\HD$.

\smallskip\noindent\emph{Diffusion.} Writing $p=\tp+p_D$, we obtain
\[
I_1=-\epsilon_+\int_\Om\Sigma_\zeta''(\tp)\,|\nabla\tp|^2\dx-\epsilon_+\int_\Om\nabla p_D\cdot\nabla\Sigma_\zeta'(\tp)\dx .
\]
The first integral is nonnegative, since $\Sigma_\zeta''\ge0$. Since $\nabla p_D\in W^{1,\sigma}(\Om)^d$ has a trace on $\partial\Om$, Green's formula applies to the second integral, and it produces no boundary term, because $\Sigma_\zeta'(\tp)=0$ on $\GD$ and $\bnu\cdot\nabla p_D=0$ on $\GN$ by (A3). Hence, by $\abs{\Sigma_\zeta'}\le1$,
\[
I_1\le-\epsilon_+\int_\Om\nabla p_D\cdot\nabla\Sigma_\zeta'(\tp)\dx=\epsilon_+\int_\Om\Sigma_\zeta'(\tp)\,\Delta p_D\dx\le\epsilon_+\norm{\Delta p_D}_1 .
\]

\smallskip\noindent\emph{Drift.} Again with $p=\tp+p_D$,
\[
I_2=-\mu_+\int_\Om\nabla S_\zeta(\tp)\cdot\nabla\phi\dx-\mu_+\int_\Om p_D\nabla\phi\cdot\nabla\Sigma_\zeta'(\tp)\dx=:I_2'+I_2'' .
\]
Since $\nabla\phi=\nabla\psi+\nabla\phi_D$ and $S_\zeta(\tp)\in\HD$, \eqref{eq:weak-poisson} and the harmonicity of $\phi_D$ in the sense of (A2) (which makes the $\phi_D$-part vanish) give
\[
I_2'=-\frac{\mu_+}{\epsilon_0}\int_\Om qS_\zeta(\tp)\dx,\qquad\text{hence}\qquad \abs{I_2'}\le\frac{\mu_+\zeta}{\epsilon_0}\norm{q}_1 .
\]

Since $p_D\in L^\infty(\Om)$, $\nabla p_D\in L^3(\Om)$ and $\phi\in H^2(\Om)\subset W^{1,6}(\Om)$, we have $p_D\nabla\phi\in H^1(\Om)^d$. As $\Sigma_\zeta'(\tp)=0$ on $\GD$ and $\bnu\cdot\nabla\phi=0$ on $\GN$, Green's formula applied to $I_2''$ produces no boundary term. Hence, by \eqref{eq:poisson},
\begin{gather*}
I_2''=\mu_+\int_\Om\Sigma_\zeta'(\tp)\nabla\cdot\big(p_D\nabla\phi\big)\dx=\mu_+\int_\Om\Sigma_\zeta'(\tp)\Big(\nabla p_D\cdot\nabla\phi-\frac{p_Dq}{\epsilon_0}\Big)\dx,\\
\abs{I_2''}\le\mu_+\norm{\nabla p_D}_2\norm{\nabla\phi}_2+\frac{\mu_+}{\epsilon_0}\norm{p_D}_2\norm{q}_2 .
\end{gather*}

\smallskip\noindent\emph{Convection.} With $p=\tp+p_D$, 
$$I_3=\int_\Om\bv\cdot\nabla S_\zeta(\tp)\dx+\int_\Om p_D\bv\cdot\nabla\Sigma_\zeta'(\tp)\dx.$$ 
Since $\nabla\cdot \bv=0$ in $\Om$ and $\bnu\cdot \bv=0$ on $\partial\Om$, the divergence theorem shows  that
\[\int_\Om\bv\cdot\nabla S_\zeta(\tp)\dx=0\quad\text{and}\quad 
I_3=-\int_\Om\Sigma_\zeta'(\tp)\bv\cdot\nabla p_D\dx,\quad\text{hence}\quad
\abs{I_3}\le\norm{\bv}_\infty\norm{\nabla p_D}_1 .
\]

\smallskip\noindent\emph{Source.} By \eqref{eq:Gamma-bound}, $n\ge0$ and $\abs{\Sigma_\zeta'}\le1$, we have $$\abs{I_4}\le\int_\Om\abs{F}\dx\le\mu_-\gamma_*\int_\Om n\abs{\nabla\phi}\dx.$$

Integrating \eqref{eq:L1-start} over $(t_1,t_2)$ and collecting the bounds on $I_1,\dots,I_4$, we obtain
\begin{multline*}
\int_\Om\Sigma_\zeta(\tp(t_2))\dx\le\int_\Om\Sigma_\zeta(\tp(t_1))\dx+\int_{t_1}^{t_2}\Big[\epsilon_+\norm{\Delta p_D}_1+\norm{\bv}_\infty\norm{\nabla p_D}_1+\mu_+\norm{\nabla p_D}_2\norm{\nabla\phi}_2\\
+\frac{\mu_+}{\epsilon_0}\norm{p_D}_2\norm{q}_2+\frac{\mu_+\zeta}{\epsilon_0}\norm{q}_1+\mu_-\gamma_*\int_\Om n\abs{\nabla\phi}\dx\Big]\dt .
\end{multline*}
Now let $\zeta\rightarrow0$. On the left, $\Sigma_\zeta(u)=\sqrt{u^2+\zeta^2}-\zeta$ increases to $\abs{u}$, so the integral converges to $\int_\Om\abs{\tp(t_2)}\dx$ by monotone convergence. On the right, $\Sigma_\zeta(u)\le\abs{u}$, and the term with $\zeta$ tends to zero because $\norm{q}_1\le\abs{\Om}^{1/2}\norm{q}_2$ and $q\in L^\infty((0,T),L^2(\Om))$. Hence, for all $0\le t_1\le t_2\le T$,
\begin{multline}\label{eq:L1-ion}
\int_\Om\abs{\tp(t_2)}\dx\le\int_\Om\abs{\tp(t_1)}\dx+\int_{t_1}^{t_2}\Big[\epsilon_+\norm{\Delta p_D}_1+\norm{\bv}_\infty\norm{\nabla p_D}_1+\mu_+\norm{\nabla p_D}_2\norm{\nabla\phi}_2\\
+\frac{\mu_+}{\epsilon_0}\norm{p_D}_2\norm{q}_2+\mu_-\gamma_*\int_\Om n\abs{\nabla\phi}\dx\Big]\dt .
\end{multline}

\smallskip\noindent\emph{The electron equation.} The same computation with $f=\Sigma_\zeta'(\tn(t))$ in \eqref{eq:weak-n} differs only in the sign of the drift term. The bounds are identical, and give \eqref{eq:L1-ion} with $(\tp,\epsilon_+,\mu_+,p_D)$ replaced by $(\tn,\epsilon_-,\mu_-,n_D)$, the source term $\mu_-\gamma_*\int_\Om n\abs{\nabla\phi}\dx$ being unchanged.

\smallskip\noindent\emph{Conclusion.} Adding \eqref{eq:L1-ion} and its counterpart for $\tn$, and collecting the coefficients, we obtain, for all $0\le t_1\le t_2\le T$,
\[\begin{split}
\int_\Om\big(\abs{\tp(t_2)}+\abs{\tn(t_2)}\big)\dx\le&\int_\Om\big(\abs{\tp(t_1)}+\abs{\tn(t_1)}\big)\dx\\
+&\int_{t_1}^{t_2}\Big(k_0+k_1\norm{\nabla\phi}_2+k_2\norm{q}_2+2\mu_-\gamma_*\int_\Om n\abs{\nabla\phi}\dx\Big)\dt\end{split}
\]
with the nonnegative constants
\begin{gather*}
k_0:=\epsilon_+\norm{\Delta p_D}_1+\epsilon_-\norm{\Delta n_D}_1+\norm{\bv}_\infty\big(\norm{\nabla p_D}_1+\norm{\nabla n_D}_1\big),\\
k_1:=\mu_+\norm{\nabla p_D}_2+\mu_-\norm{\nabla n_D}_2,\qquad
k_2:=\frac{\mu_+\norm{p_D}_2+\mu_-\norm{n_D}_2}{\epsilon_0},
\end{gather*}
and we set
\begin{equation}\label{eq:C3}
C_3:=k_0+k_1+k_2 .
\end{equation}
Since $k_0+k_1\norm{\nabla\phi}_2+k_2\norm{q}_2\le C_3\big(1+\norm{\nabla\phi}_2+\norm{q}_2\big)$, this is \eqref{eq:L1}.
\end{proof}

\subsection{Closure of the estimates}\label{sec:proofB}

\begin{proposition}[Closure]\label{prop:closure}
Let Assumption~\ref{ass:main}, \textup{(H1)} and \textup{(H2)} hold and let $(\phi,p,n)$ be a weak solution on $(0,T)$. Then there is a constant $C>0$, depending only on $\mathcal D$, such that, for all $t\in[0,T]$,
\begin{equation}\label{eq:closure}
\ee(t)+m(t)+\frac12\int_0^t\Jou\dr+\frac{\epsilon}{2\epsilon_0}\int_0^t\norm{q}_2^2\dr\le C(1+t)e^{Ct}\big(1+\ee(0)+m(0)\big).
\end{equation}
\end{proposition}

\begin{proof}
We first bound the source term of \eqref{eq:L1} by the Joule heating and the number of carriers. Since $\nabla\phi=\nabla\psi+\nabla\phi_D$, and since $p,n\ge0$ gives $\int_\Om n\dx\le\int_\Om s\dx$ and $\mu_-n\le g$, the Cauchy--Schwarz inequality with respect to the measure $n\dx$ and (H2) give
\begin{equation}\label{eq:source-Joule}
\begin{aligned}
\int_\Om n\abs{\nabla\phi}\dx&\le\Big(\int_\Om n\dx\Big)^{1/2}\Big(\int_\Om n\abs{\nabla\psi}^2\dx\Big)^{1/2}+\norm{\nabla\phi_D}_\infty\int_\Om n\dx\\
&\le\mu_-^{-1/2}\Jou^{1/2}\Big(\int_\Om s\dx\Big)^{1/2}+\norm{\nabla\phi_D}_\infty\int_\Om s\dx .
\end{aligned}
\end{equation}
Set $Z(t):=\ee(t)+m(t)$, which is continuous on $[0,T]$. By \eqref{eq:mass}, 
$$\int_\Om s\dx\le m+C\le Z+C.$$ 
Also 
$$\norm{\nabla\phi}_2\le\norm{\nabla\phi_D}_2+\norm{\nabla\psi}_2\le C(1+\ee).$$

\smallskip\noindent\emph{Energy.} Integrating \eqref{eq:energy-ineq} over $(0,t)$, we obtain
\begin{equation}\label{eq:E-int}
\ee(t)+\frac12\int_0^t\Jou\dr+\frac{\epsilon}{2\epsilon_0}\int_0^t\norm{q}_2^2\dr\le\ee(0)+C\int_0^t\big(1+Z\big)\dr .
\end{equation}

\smallskip\noindent\emph{Mass.} We take $t_1=0$ in \eqref{eq:L1} of Lemma~\ref{lem:L1} and insert \eqref{eq:source-Joule}, in which we bound the product $\Jou^{1/2}(\int_\Om s\dx)^{1/2}$ by Young's inequality,
\[
\Jou^{1/2}\Big(\int_\Om s\dx\Big)^{1/2}\le\frac12\Jou+\frac12\int_\Om s\dx .
\]
With $\int_\Om s\dx\le Z+C$ and $\norm{\nabla\phi}_2\le C(1+\ee)$, this gives
\[
m(t)\le m(0)+C\int_0^t\big(1+Z\big)\dr+C\int_0^t\norm{q}_2\dr+\mu_-^{1/2}\gamma_*\int_0^t\Jou\dr ,
\]
the coefficient of $\Jou$ being $2\mu_-\gamma_*\cdot\mu_-^{-1/2}\cdot\frac12=\mu_-^{1/2}\gamma_*$.

By \eqref{eq:E-int}, 
$$\int_0^t\Jou\dr\le2\ee(0)+C\int_0^t(1+Z)\dr,$$ 
while by the Cauchy--Schwarz inequality, \eqref{eq:E-int} and $\sqrt{ab}\le\frac12(a+b)$,
\[
\int_0^t\norm{q}_2\dr\le t^{1/2}\Big(\int_0^t\norm{q}_2^2\dr\Big)^{1/2}\le\frac t2+\frac{\epsilon_0}{\epsilon}\Big(\ee(0)+C\int_0^t(1+Z)\dr\Big).
\]
Hence
\begin{equation}\label{eq:m-int}
m(t)\le m(0)+C\big(1+t\big)\big(1+\ee(0)\big)+C\int_0^tZ\dr .
\end{equation}

\smallskip\noindent\emph{Gronwall.} Adding \eqref{eq:E-int} and \eqref{eq:m-int} and dropping the two nonnegative integrals on the left of \eqref{eq:E-int} gives, for all $t\in[0,T]$,
\[
Z(t)\le C(1+t)\big(1+Z(0)\big)+C\int_0^tZ\dr ,\qquad Z(0)=\ee(0)+m(0).
\]
Since $Z$ is continuous and $C(1+t)(1+Z(0))$ is nondecreasing in $t$, Gronwall's inequality yields
\begin{equation}\label{est_enegy+mass}
Z(t)\le C(1+t)\big(1+Z(0)\big)e^{Ct},
\end{equation}
which bounds $\ee(t)+m(t)$ by the right-hand side of \eqref{eq:closure}. 
The two remaining terms on the left of \eqref{eq:closure}, $\frac12\int_0^t\Jou\dr$ and $\frac{\epsilon}{2\epsilon_0}\int_0^t\norm{q}_2^2\dr$, also stand on the left of \eqref{eq:E-int}. Dropping $\ee(t)\ge0$ on the left of \eqref{eq:E-int} and inserting the bound \eqref{est_enegy+mass} on $Z$ into the integral $C\int_0^t(1+Z)\dr$ on the right of \eqref{eq:E-int}, we obtain, since $1+r\le1+t$ for $0\le r\le t$,
\begin{align*}
\int_0^t\big(1+Z\big)\dr&\le t+C\big(1+Z(0)\big)\int_0^t(1+r)e^{Cr}\dr\\
&\le t+C\big(1+Z(0)\big)(1+t)\int_0^te^{Cr}\dr\le C(1+t)e^{Ct}\big(1+Z(0)\big),
\end{align*}
and hence
\[
\frac12\int_0^t\Jou\dr+\frac{\epsilon}{2\epsilon_0}\int_0^t\norm{q}_2^2\dr\le\ee(0)+C\int_0^t\big(1+Z\big)\dr\le C(1+t)e^{Ct}\big(1+\ee(0)+m(0)\big).
\]
Together with the bound \eqref{est_enegy+mass} on $\ee(t)+m(t)$, this proves \eqref{eq:closure}.
\end{proof}

\begin{proof}[Proof of Theorem~\ref{thm:B}]
Let $(\phi,p,n)$ be the maximal weak solution and suppose $T_{\max}<\infty$. The numbers $\ee(0)=\norm{p_0-n_0}_{\HDd}^2/(2\epsilon_0)$ and $m(0)=\norm{p_0-p_D}_1+\norm{n_0-n_D}_1$ are finite and depend only on $p_0,n_0$ and the data. For every $T'<T_{\max}$, Proposition~\ref{prop:closure}, applied to the restriction of the solution to $(0,T')$, gives 
$$\sup_{[0,T']}\ee\le C(1+T_{\max})e^{CT_{\max}}(1+\ee(0)+m(0))=:K_{\ee}.$$ 
Hence $\ee\le K_{\ee}$ on $[0,T_{\max})$. Since $\norm{\nabla\phi}_2^{2\kappa}\le C(1+\ee)^\kappa\le C(1+K_{\ee})^\kappa$, we obtain $$\int_0^{T_{\max}}\norm{\nabla\phi}_2^{2\kappa}\dt<\infty,$$ which contradicts Theorem~\ref{thm:A}. Therefore $T_{\max}=\infty$.

The bounds \eqref{eq:growth}, with the dependence of $C$ stated there, are established in Corollary~\ref{cor:growth} below.
\end{proof}

\begin{corollary}[Growth bounds]\label{cor:growth}
Under the assumptions of Theorem~\ref{thm:B} there is $C$, depending only on $\mathcal D$, such that for all $t\ge0$
\begin{subequations}\label{bd4thm12}
\begin{equation}\label{1stbd}\ee(t)+\int_\Om\big(p(t)+n(t)\big)\dx\le C\big(1+\ee(0)+m(0)\big)(1+t)e^{Ct},\end{equation}\begin{equation}\label{2ndbd}
\norm{p(t)}_2^2+\norm{n(t)}_2^2+\int_0^t\big(\norm{\nabla p}_2^2+\norm{\nabla n}_2^2\big)\dr\le C\big(1+\Y(0)\big)\exp\Big(C\big(1+\ee(0)+m(0)\big)^\kappa e^{Ct}\Big).\end{equation}
\end{subequations}
In particular, \eqref{eq:growth} holds with $N_0=1+\norm{p_0}_2^2+\norm{n_0}_2^2$ as in Theorem~\ref{thm:B}. The constant $C$ can be chosen nondecreasing in $\gamma_*=\max\{\alpha_1,\eta_0\}$ and in each of the norms of $\phi_D$, $p_D$, $n_D$ and $\bv$ listed in Theorem~\ref{thm:B}.
\end{corollary}

\begin{proof}
Since $T_{\max}=\infty$, Proposition~\ref{prop:closure} applies on $(0,T)$ for every $T>0$. The estimate \eqref{1stbd} is \eqref{eq:closure} together with \eqref{eq:mass}. For the estimate \eqref{2ndbd}, we bound the coefficient $a$ of \eqref{eq:coefficient}. Since 
\[\norm{\nabla\phi}_2\le\norm{\nabla\phi_D}_2+(2\ee/\epsilon_0)^{1/2},\] we have \[\norm{\nabla\phi}_2^{2\kappa}\le C(1+\ee)^\kappa,\] and the estimate \eqref{1stbd} gives \[1+\ee(r)\le C\big(1+\ee(0)+m(0)\big)(1+r)e^{Cr},\] hence
\begin{equation}\label{int_a}
\int_0^ta(r)\dr\le C_1t+C\big(1+\ee(0)+m(0)\big)^\kappa\int_0^t(1+r)^\kappa e^{C\kappa r}\dr\le C\big(1+\ee(0)+m(0)\big)^\kappa e^{Ct},
\end{equation}
where we used $(1+r)^\kappa e^{C\kappa r}\le e^{(C+1)\kappa r}$ and $C_1t\le C_1e^{t}$.

Inserting \eqref{int_a} into \eqref{eq:field-gronwall} gives
\[
1+\Y(t)+\frac\theta2\int_0^t\Big(\norm{\nabla\tp}_2^2+\norm{\nabla\tn}_2^2\Big)\dr
\le\big(1+\Y(0)\big)\exp\Big(C\big(1+\ee(0)+m(0)\big)^\kappa e^{Ct}\Big)=:Q(t).
\]
Since $\nabla p=\nabla\tp+\nabla p_D$ and $\nabla n=\nabla\tn+\nabla n_D$, we have $$\norm{\nabla p}_2^2+\norm{\nabla n}_2^2\le2\big(\norm{\nabla\tp}_2^2+\norm{\nabla\tn}_2^2\big)+C,$$ and $\norm{p}_2^2+\norm{n}_2^2\le C(1+\Y)$ by \eqref{eq:Y}. Hence
\[
\norm{p(t)}_2^2+\norm{n(t)}_2^2+\int_0^t\Big(\norm{\nabla p}_2^2+\norm{\nabla n}_2^2\Big)\dr\le C\,Q(t)+Ct\le C\,Q(t),
\]
where, enlarging $C$ if necessary, we have assumed $C\ge1$, so that $t\le e^{Ct}\le Q(t)$.

Finally, by \eqref{eq:weak-poisson}, \eqref{energy=Poisson} and the Poincar\'e inequality $\norm{f}_{\HDd}\le C_P\norm{f}_2$,
\begin{gather*}
\ee(0)=\tfrac1{2\epsilon_0}\norm{p_0-n_0}_{\HDd}^2\le C\big(\norm{p_0}_2^2+\norm{n_0}_2^2\big),\\
m(0)\le\abs\Om^{1/2}\big(\norm{p_0}_2+\norm{n_0}_2\big)+\norm{p_D}_1+\norm{n_D}_1,
\end{gather*}
and $\Y(0)\le C\big(1+\norm{p_0}_2^2+\norm{n_0}_2^2\big)$, which gives \eqref{eq:growth}.

It remains to verify the dependence of $C$ asserted in Theorem~\ref{thm:B} and the monotonicity asserted in the corollary. By the convention following \eqref{eq:data}, $C$ depends only on $\mathcal D$; in particular, it is independent of the initial data and of time. Under \textup{(H1)}--\textup{(H2)}, the entries of $\mathcal D$ are the arguments of $C$ in Theorem~\ref{thm:B}, except that $\mathcal D$ contains the three constants $\alpha_1,\alpha_2,\eta_0$ where Theorem~\ref{thm:B} has only $\gamma_*=\max\{\alpha_1,\eta_0\}$. The constants $\alpha_1,\alpha_2,\eta_0$ enter the estimates of Sections~\ref{sec:L2} and~\ref{sec:global} that determine $C$ only through \eqref{eq:Gamma-bound}, that is, through $\gamma_*$; in particular, $C$ does not depend on $\alpha_2$. Moreover, $\gamma_*$ and the norms of $\phi_D$, $p_D$, $n_D$ and $\bv$ that appear in Theorem~\ref{thm:B} enter these estimates only as upper bounds, so each estimate remains valid when any of them is replaced by a larger number. Hence $C$ can be chosen nondecreasing in each of them. Remark~\ref{rem:constants} below records where each of these quantities enters.
\end{proof}
The second bound in \eqref{eq:growth} is double-exponential in $t$ because the exponential bound of Proposition~\ref{prop:closure} on $\ee$ enters the exponent of \eqref{eq:field-gronwall}. In the simulations of Section~\ref{sec:numerics} the number of carriers grows by orders of magnitude and then levels off or decreases (Figure~\ref{fig:voltages}(a)). Sharpening the growth rates in
\eqref{eq:growth} is an open question which we do not pursue here.

\begin{remark}[What the constant depends on]\label{rem:constants}
The constant $C$ of \eqref{eq:growth} is assembled in Proposition~\ref{prop:closure}, \eqref{eq:field-gronwall} and Corollary~\ref{cor:growth}. Its ingredients are of the following four kinds; a last paragraph records how the initial data and the time enter \eqref{eq:growth}.

\emph{Geometry.} The dimension $d$ and the domain $\Om$, the latter only through $\abs\Om$, the constant $C_{GN}$ of Lemma~\ref{lem:aux}(i), the Poincar\'e constant of $\HD$ and the constants of the Sobolev embeddings of Section~\ref{sec:setting}, among them that in \eqref{eq:v-infty}. In the estimates that determine $C$ (Lemmas~\ref{lem:L2identity}, \ref{lem:data} and~\ref{lem:source}(ii), Proposition~\ref{prop:field-controlled} and Section~\ref{sec:global}), the $H^2$-regularity of Lemma~\ref{lem:aux}(ii), and with it the angle condition in (A1), is used only qualitatively, to justify integrations by parts.

\emph{Physical constants.} The diffusivity $\epsilon$, the mobilities $\mu_\pm$ and the permittivity $\epsilon_0$. The constants $\alpha_1,\alpha_2,\eta_0$ enter \emph{only} through $\gamma_*=\max\{\alpha_1,\eta_0\}$, since the estimates of Sections~\ref{sec:L2} and~\ref{sec:global} use no property of $\Gamma$ beyond \eqref{eq:Gamma-bound}. Hence $C$ does not depend on $\alpha_2$: for fixed $\alpha_1$ and $\eta_0$, the bounds \eqref{eq:growth} hold with one constant for all $\alpha_2>0$. When $\alpha_1>\eta_0$, the net Townsend coefficient $\alpha_1e^{-\alpha_2/E}-\eta_0$ is negative for $E$ below the threshold field $\alpha_2/\log(\alpha_1/\eta_0)$ and positive above it. The threshold is proportional to $\alpha_2$, so the bounds hold with one constant for every value of the threshold field. Nor do the bounds depend on the sign of $\Gamma$: they use only $|\Gamma(E)|\le\mu_-\gamma_*E$, so they are the same for fields below the threshold, where attachment prevails, and above it, where ionization prevails. By contrast, $\alpha_2$ matters in Proposition~\ref{prop:scaling} and Section~\ref{ssec:touchdown}, and the threshold field in Section~\ref{sec:numerics}.

\emph{Boundary data.} The boundary data enter through $\norm{\phi_D}_{H^2}$, through $\norm{\nabla\phi_D}_{L^\infty}$ under (H2), and through $\norm{p_D}_{W^{2,\sigma}}$ and $\norm{n_D}_{W^{2,\sigma}}$. Of the carrier data, the estimates of Sections~\ref{sec:L2} and~\ref{sec:global} use quantitatively only $\norm{p_D}_{L^\infty}$, $\norm{\nabla p_D}_{L^1}$, $\norm{\nabla p_D}_{L^2}$ and $\norm{\Delta p_D}_{L^1}$, their counterparts for $n_D$, and the $L^1$- and $L^2$-norms of $p_D,n_D$, which are bounded by their $L^\infty$-norms and $\abs\Om$; the remaining properties in (A3), such as $\nabla p_D\in L^3(\Om)$ and the Neumann compatibility on $\GN$, serve only to justify integrations by parts. The constants $R_a,R_b$ of (A3) do not enter: they belong to the hypotheses of the local theory \cite[Assumption~1.5]{GHL26} and are not used in Sections~\ref{sec:L2} and~\ref{sec:global}.

\emph{Gas velocity.} The velocity enters only through $\norm{\bv}_\infty$ of \eqref{eq:v-infty}, hence through $\norm{\bv}_{L^\infty((0,\infty),H^2)}$. Since (A4) bounds $\bv$ on the whole half-line, this norm does not depend on the time interval, so the length of the interval does not enter $C$ through the velocity. In the estimates, the conditions $\nabla\cdot\bv=0$ and $\bnu\cdot\bv=0$ in (A4) serve only to handle the convection terms in Lemmas~\ref{lem:L2identity} and~\ref{lem:L1}, where the parts carrying $\tp$ and $\tn$ vanish.

\emph{Initial data and time.} Neither enters $C$. The initial data enter \eqref{eq:growth} only through $N_0$, and the time only through the factors $(1+t)e^{Ct}$ and $\exp(CN_0^\kappa e^{Ct})$, which come from the Gronwall arguments of Proposition~\ref{prop:closure} and of \eqref{eq:field-gronwall}. The right-hand sides of \eqref{eq:growth} are nondecreasing in $C$, and $C$ is nondecreasing in $\gamma_*$ and in the norms of the boundary data and of the velocity (Corollary~\ref{cor:growth}). Hence, for fixed $d,\Om,\epsilon,\mu_\pm,\epsilon_0$, one constant $C$ serves for all boundary data, velocity fields and ionization constants for which $\gamma_*$ and these norms stay below given bounds, and for all initial data. No smallness is required: \eqref{eq:growth} holds for every voltage admitted by (A2) and (H2) and every initial state admitted by (A3).
\end{remark}

\begin{remark}[The role of \textup{(H1)}]\label{rem:H1}
Hypothesis (H1) is used at a single point: in the proof of Corollary~\ref{cor:energy-ineq}, it makes the first term on the right of \eqref{eq:energy-identity}, the ambipolar term of Remark~\ref{rem:interpretation}, vanish. Theorem~\ref{thm:A}, Proposition~\ref{prop:field-controlled} and Lemmas~\ref{lem:energy-identity} and~\ref{lem:L1} hold for arbitrary $\epsilon_\pm>0$. Assume (H2) but not (H1). Then the ambipolar term need not vanish, and we denote it by
\begin{equation}\label{eq:ambipolar}
\mathcal A:=-\frac{\epsilon_+-\epsilon_-}{2\epsilon_0}\int_\Om\big(p^2-n^2\big)\dx=\frac{\epsilon_--\epsilon_+}{2\epsilon_0}\int_\Om q\,s\dx .
\end{equation}
The proof of Corollary~\ref{cor:energy-ineq} then gives, for every $\vartheta\in(0,1)$, a constant $C_\vartheta$ depending only on $\mathcal D$ and $\vartheta$ such that for a.e.\ $t\in(0,T)$
\begin{equation}\label{eq:energy-ineq-general}
\frac{d\ee}{dt}+\vartheta\Jou+\vartheta\,\frac{\epsilon_++\epsilon_-}{2\epsilon_0}\norm{q}_2^2\le C_\vartheta\Big(1+\ee+\int_\Om s\dx\Big)+\mathcal A .
\end{equation}
Indeed, Young's inequality with the weight $1-\vartheta$ in place of $\frac12$ bounds the battery term by $(1-\vartheta)\Jou+C_\vartheta\int_\Om s\dx$. Since the coefficient of $\norm{q}_2^2$ on the left of \eqref{eq:energy-identity} is $\frac{\epsilon_++\epsilon_-}{2\epsilon_0}$, Young's inequality in the data term $\frac1{\epsilon_0}\int_\Om(\epsilon_+p_D-\epsilon_-n_D)\,q\dx$ and in the convection term, with the weight $(1-\vartheta)\frac{\epsilon_++\epsilon_-}{4\epsilon_0}$ in place of $\frac{\epsilon}{4\epsilon_0}$, leaves $\vartheta\frac{\epsilon_++\epsilon_-}{2\epsilon_0}\norm{q}_2^2$ on the left and adds $C_\vartheta(1+\ee)$ to the right. The data term $-\int_\Om\nabla(\epsilon_+p_D-\epsilon_-n_D)\cdot\nabla\psi\dx$ is bounded by $C_\vartheta(1+\ee)$ as in the proof of Corollary~\ref{cor:energy-ineq}. For $\vartheta=\frac12$ and under (H1), \eqref{eq:energy-ineq-general} is \eqref{eq:energy-ineq}.

The proof of Proposition~\ref{prop:closure} therefore extends to $\epsilon_+\neq\epsilon_-$ as soon as $\int_0^t\mathcal A\dr$ is controlled by the Joule heating and the dissipative term. Precisely, we ask for constants $\delta<\frac{\epsilon_++\epsilon_-}{2\epsilon_0}$, $\beta<1$ and $C$, depending only on $\mathcal D$, such that every weak solution on every interval $(0,T)$ satisfies, for $0\le t\le T$,
\begin{equation}\label{eq:needed}
\int_0^t\mathcal A\dr\le\delta\int_0^t\norm{q}_2^2\dr+\beta\int_0^t\Jou\dr+C(1+t)\big(1+\ee(0)+m(0)\big)+C\int_0^t\big(\ee+m\big)\dr .
\end{equation}
Under (H1), $\mathcal A=0$, so \eqref{eq:needed} holds with $\delta=\beta=0$, and Proposition~\ref{prop:conditional} below reduces to Theorem~\ref{thm:B}. Whether \eqref{eq:needed} holds for $\epsilon_+\neq\epsilon_-$ is open: Section~\ref{ssec:num-uneq} evaluates its terms on a computed solution, and Section~\ref{ssec:unequal} discusses why it has resisted proof.
\end{remark}

\begin{proposition}[Conditional global existence]\label{prop:conditional}
Let Assumption~\ref{ass:main} and \textup{(H2)} hold, and suppose that \eqref{eq:needed} holds with constants $\delta$, $\beta$ and $C$ as in Remark~\ref{rem:H1}. Then $T_{\max}=\infty$, and the bounds \eqref{eq:growth} hold with a constant depending on $\mathcal D$, $\delta$ and $\beta$.
\end{proposition}
\begin{proof}
Choose $\vartheta\in(\max\{\beta,0\},1)$ with $\vartheta\frac{\epsilon_++\epsilon_-}{2\epsilon_0}>\delta$.
Integrating \eqref{eq:energy-ineq-general} over $(0,t)$ and inserting \eqref{eq:needed} gives \eqref{eq:E-int}
with the coefficients $\vartheta-\beta$ and $\vartheta\frac{\epsilon_++\epsilon_-}{2\epsilon_0}-\delta$
in place of $\frac12$ and $\frac{\epsilon}{2\epsilon_0}$, and with the additional term
$C(1+t)(1+\mathcal E(0)+m(0))$ on the right. The Mass and Gronwall steps of the proof of Proposition~\ref{prop:closure} use only Lemma~\ref{lem:L1}, \eqref{eq:source-Joule} and \eqref{eq:E-int}, the last in the form just obtained; none of them requires (H1), and they go through with modified constants, as do the proofs of Theorem~\ref{thm:B} and Corollary~\ref{cor:growth}.
\end{proof}

\section{Numerical experiments}\label{sec:numerics}

In this section we report two numerical experiments on the questions that Theorems~\ref{thm:A} and~\ref{thm:B} leave open, for a discharge in the MEMS capacitor of \cite{GHL26} computed with the finite volume scheme of Appendix~\ref{app:scheme}: we examine how the discharge changes with the voltage, in connection with Problem~\ref{conj:threshold}, and we evaluate the terms of the bound \eqref{eq:needed}, to which Proposition~\ref{prop:conditional} reduces global existence for unequal diffusivities, in a simulation with the diffusivities of air.

\subsection{Setting}\label{ssec:num-setting}
\begin{figure}[htbp]
\centering
\includegraphics[width=0.6\textwidth]{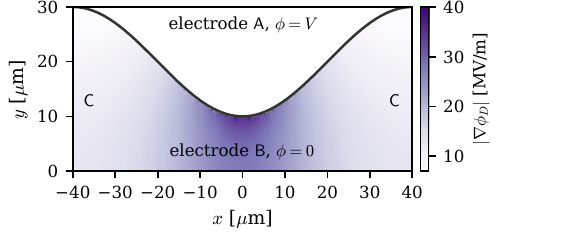}
\caption{The computational domain \eqref{eq:num-domain} at a stage before touch-down, with the modulus of the vacuum field $\nabla\phi_D$ at $V=300$ V: the anode $y=w(x)$ is the deformable electrode $\mathsf A$ at potential $V$, the cathode $y=0$ the fixed electrode $\mathsf B$ at potential $0$, and the Neumann boundary $\GN$ at $|x|=r$ the artificial boundary $\mathsf C$. The field is largest at the bottom of the dip.}\label{fig:geometry}
\end{figure}

The MEMS capacitor of \cite[Figure~1]{GHL26} consists of a deformable electrode $\mathsf A$ at potential $V$, a fixed electrode $\mathsf B$ at potential $0$ and an artificial boundary $\mathsf C$ at $|x|=r$, on which Neumann conditions are imposed. We take this geometry as the computational domain (Figure~\ref{fig:geometry}), at a stage before touch-down at which $\mathsf A$ has been pulled towards $\mathsf B$, so that the gap is narrowest at $x=0$; as in \cite{GHL26}, the position of the electrodes is prescribed. Thus
\begin{equation}\label{eq:num-domain}
\Om=\{(x,y):\ |x|<r,\ 0<y<w(x)\},\qquad w(x)=w_{\min}+\frac{w_{\max}-w_{\min}}{2}\Big(1-\cos\frac{\pi x}{r}\Big),
\end{equation}
with the cathode $\mathsf B=\{y=0\}$, the anode $\mathsf A=\{y=w(x)\}$ and $\GN=\{|x|=r\}$ in the role of $\mathsf C$. The profile $w$ attains its minimum $w_{\min}$ at $x=0$ and its maximum $w_{\max}$ on $\mathsf C$, where $w'(\pm r)=0$, so that the electrodes meet $\GN$ at right angles, which (A1) admits; the discharge is controlled mainly by $w_{\min}$ and the width of the dip, since the field is close to $V/w(x)$ where $w$ varies slowly. Hypothesis (H2) holds: as $w$ is smooth with $w'(\pm r)=0$, the even reflection of $\phi_D$ across $|x|=r$ is harmonic in a smooth periodic strip with constant data on each electrode, so $\nabla\phi_D$ is bounded (Section~\ref{ssec:touchdown}). We compute $\phi_D$ on the same grid as the solution (Figure~\ref{fig:geometry}); on the grid, $\max\abs{\nabla\phi_D}=1.35\,V/w_{\min}$ (Appendix~\ref{app:scheme}). The gas is at rest, $\bv=0$, since in a MEMS gap the gas moves at metres per second while the carriers drift at $10^3$ to $10^6$ metres per second; it is air at atmospheric pressure, with the Townsend constants and the mobilities of Table~\ref{tab:params}, which also lists the geometry. The diffusivities of air do not satisfy (H1); Section~\ref{ssec:num-voltage} uses the common value $\epsilon=\epsilon_+=\epsilon_-$ of Table~\ref{tab:params}, and Section~\ref{ssec:num-uneq} repeats the simulation at $300$ V with the diffusivities of air. The electrode data are the constants $p_D=n_D=n_b$ of \cite{GHL26}, which satisfy (A3) with $R_a=R_b=n_b$ and carry no charge, $q_D=0$; the localization of the discharge is thus due to the geometry alone. The initial data are $p_0=n_0=n_b+\varsigma$ with the seed $\varsigma(x,y)=\varsigma_0\,e^{-(x^2+(y-y_s)^2)/(2b^2)}\,4y\,(w_{\min}-y)_+/w_{\min}^2$, a Gaussian with the prefactor $\varsigma_0=10^{-2}$ C/m$^3$ and the width $b=1$ $\mu$m, centred on the axis of the dip at the height $y_s=3$ $\mu$m above the cathode and cut off by the last factor so that it vanishes on both electrodes. The seed is Lipschitz continuous but not $C^1$, so the initial data satisfy \eqref{eq:restart-data} but not (A3); by Proposition~\ref{prop:restart}, the maximal weak solution of Section~\ref{ssec:local} exists for such data, and Theorems~\ref{thm:A} and~\ref{thm:B} hold for it (Section~\ref{ssec:local}).

\begin{table}[htbp]
\centering
\small
\begin{tabular}{llll}
\toprule
Quantity & Symbol & Value & Unit\\
\midrule
gap width at the dip, $w(0)$ & $w_{\min}$ & $10$ & $\mu$m\\
gap width on $\mathsf C$, $w(\pm r)$ & $w_{\max}$ & $30$ & $\mu$m\\
half-length of the domain & $r$ & $40$ & $\mu$m\\
applied voltage & $V$ & $150$, $300$, $450$ & V\\
electron mobility & $\mu_-$ & $4\cdot10^{-2}$ & m$^2$/(V\,s)\\
ion mobility & $\mu_+$ & $2\cdot10^{-4}$ & m$^2$/(V\,s)\\
diffusivities, simulations with (H1) & $\epsilon_+=\epsilon_-$ & $10^{-4}$ & m$^2$/s\\
diffusivities, simulation of Section~\ref{ssec:num-uneq} & $\epsilon_+$, $\epsilon_-$ & $5\cdot10^{-6}$, $8\cdot10^{-2}$ & m$^2$/s\\
Townsend coefficient & $\alpha_1$ & $1.14\cdot10^{6}$ & 1/m\\
Townsend field constant & $\alpha_2$ & $2.774\cdot10^{7}$ & V/m\\
attachment coefficient & $\eta_0$ & $5\cdot10^{3}$ & 1/m\\
electrode data and background density & $n_b$ & $10^{-3}$ & C/m$^3$\\
\bottomrule
\end{tabular}
\caption{Physical parameters of the simulations. The Townsend constants are $\alpha_1=A_{\mathrm T}P$ and $\alpha_2=B_{\mathrm T}P$ with the values $A_{\mathrm T}=15$~1/(cm\,Torr), $B_{\mathrm T}=365$~V/(cm\,Torr) of air \cite{Rai91} at the pressure $P=760$ Torr; the mobilities are those of electrons and of molecular ions in air at atmospheric pressure \cite{Rai91}; the diffusivities of Section~\ref{ssec:num-uneq} are, up to rounding, the values $\epsilon_\pm=\mu_\pm k_BT_\pm/e$ of the Einstein relation with the gas temperature $300$ K for the ions and an electron temperature of $2$ eV, and the attachment coefficient is a representative value for air at atmospheric pressure. The densities are charge densities, $1$ C/m$^3$ corresponding to $6.2\cdot10^{12}$ carriers per cm$^3$.}\label{tab:params}
\end{table}

In all three simulations the field exceeds the Townsend threshold in the dip: the threshold field $\alpha_2/\log(\alpha_1/\eta_0)=5.1$ MV/m lies below the mean fields $V/w_{\min}=15$, $30$ and $45$ MV/m there, and below the field $V/w_{\max}$ on the flanks at $300$ and $450$ V but not at $150$ V. The ion transit time across the dip, $w_{\min}^2/(\mu_+V)$, is between $1.1$ and $3.3$ ns; the simulations cover $10$ ns, or $5$ ns at $450$ V, and the scheme, its time step and their accuracy are described in Appendix~\ref{app:scheme}.

\subsection{Dependence on the voltage}\label{ssec:num-voltage}
\begin{figure}[htbp]
\centering
\includegraphics[width=\textwidth]{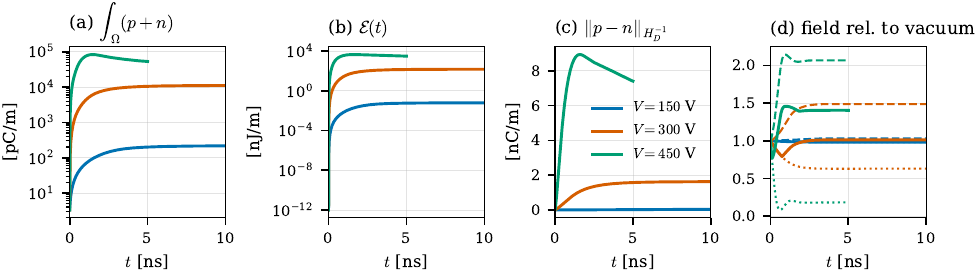}
\caption{Dependence on the voltage, for $V=150$, $300$ and $450$ V: (a) number of carriers, (b) electrostatic energy, (c) $\HDd$-norm of the charge, and (d) field strengths relative to the vacuum field, the maximal field relative to the maximal vacuum field on the grid (solid) and the fields at the cathode (dashed) and at the anode (dotted) on the axis $x=0$ relative to their vacuum values. The simulation at $450$ V is stopped at $5$ ns, after the overshoot described in Section~\ref{ssec:num-voltage}. The legend in (c) applies to all panels.}\label{fig:voltages}
\end{figure}

Figure~\ref{fig:voltages} compares the three voltages. At $150$ V the discharge is weak (Figure~\ref{fig:voltages}(a)--(c)) and the fields stay within a few percent of their vacuum values (Figure~\ref{fig:voltages}(d)); at $300$ V a positive ion cloud forms below the dip, which enhances the field at the cathode to about $1.5$ times its vacuum value and reduces it at the anode to about $0.6$ times, the beginning of a cathode fall (Figure~\ref{fig:voltages}(d)). At both voltages the number of carriers and the electrostatic energy level off within the $10$ ns of the simulations (Figure~\ref{fig:voltages}(a) and~(b)). Since at $150$ V they do so although the fields barely change, the levelling off reflects the balance between the production of the ions and their loss at the cathode and does not require a redistribution of the field. The space charge shifts this balance, at $450$ V strongly enough for the dynamics to overshoot. Within the first nanosecond the anode field on the axis falls briefly below the Townsend threshold field of Section~\ref{ssec:num-setting}, to about $8$ percent of its vacuum value, so that the net ionization there ceases; the cathode and anode fields then settle after about $2$ ns at about twice and about a fifth of their vacuum values (Figure~\ref{fig:voltages}(d)). The number of carriers, the electrostatic energy and the $\HDd$-norm of the charge pass through a maximum at about $1.5$ ns and are still decreasing when the simulation stops at $5$ ns (Figure~\ref{fig:voltages}(a)--(c)), although the fields on the axis have settled (Figure~\ref{fig:voltages}(d)). Between $150$ and $450$ V the number of carriers increases by more than two orders of magnitude and the electrostatic energy by almost five orders of magnitude (Figure~\ref{fig:voltages}(a) and~(b)), whereas the maximal field stays below $1.5$ times the maximal vacuum field (solid lines in Figure~\ref{fig:voltages}(d)) and the $\HDd$-norm of the charge below $9$ nC/m (Figure~\ref{fig:voltages}(c)), so that, since $\norm{\nabla\phi}_2^2=\norm{\nabla\phi_D}_2^2+\norm{\nabla\psi}_2^2$ by (A2), $\norm{\nabla\psi}_2=\epsilon_0^{-1}\norm{p-n}_{\HDd}$ by Lemma~\ref{lem:aux}(ii) and $\norm{\nabla\phi_D}_2=1.0$ kV at $450$ V, the norm $\norm{\nabla\phi}_2$, through which Theorem~\ref{thm:A} controls the solution, stays below $1.5$ times its vacuum value. At the end of the simulations none of the quantities in Figure~\ref{fig:voltages}(a)--(c) is growing appreciably: at $150$ and $300$ V they have levelled off, changing by less than one percent over the last nanosecond, and at $450$ V they are still decreasing when the simulation stops.

\subsection{Unequal diffusivities}\label{ssec:num-uneq}
\begin{figure}[htbp]
\centering
\includegraphics[width=\textwidth]{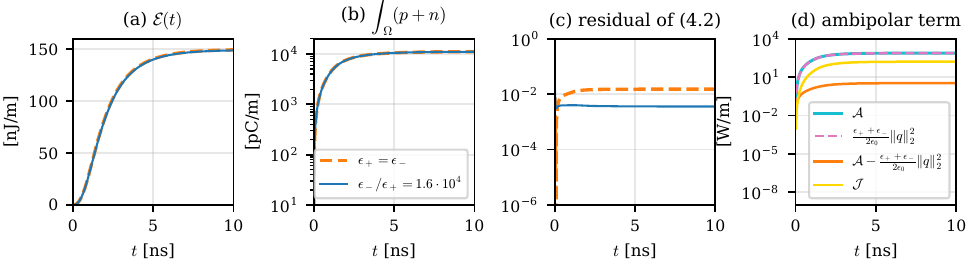}
\caption{Equal versus unequal diffusivities at $V=300$ V: (a) electrostatic energy; (b) number of carriers; (c) relative residual of the balance law \eqref{eq:energy-identity}, defined in Section~\ref{ssec:num-uneq}; with unequal diffusivities the sum of terms that normalizes it is dominated by $\mathcal A$ and the dissipative term; (d) in the simulation with unequal diffusivities, the ambipolar term $\mathcal A$ of \eqref{eq:ambipolar}, the dissipative term $\frac{\epsilon_++\epsilon_-}{2\epsilon_0}\norm q_2^2$, their difference and the Joule heating $\Jou$. The legend in (b) applies to (a)--(c); the case of equal diffusivities is dashed.}\label{fig:uneq}
\end{figure}

Without (H1) the balance law \eqref{eq:energy-identity} contains the ambipolar term $\mathcal A$ of \eqref{eq:ambipolar}, and Proposition~\ref{prop:conditional} reduces global existence to the bound \eqref{eq:needed} on that term. Figure~\ref{fig:uneq} compares the simulation at $300$ V of Section~\ref{ssec:num-voltage} with the same simulation for the diffusivities of air in Table~\ref{tab:params}, $\epsilon_-/\epsilon_+=1.6\cdot10^4$. The electrostatic energy and the number of carriers, the two quantities through which Proposition~\ref{prop:closure} bounds the solution, differ by less than two percent in the two simulations (Figure~\ref{fig:uneq}(a) and~(b)): in this simulation, (H1) makes no visible difference to what Theorem~\ref{thm:B} controls. After the initial transient the balance law holds in both simulations up to a relative residual of $1.5$ percent with equal and $0.4$ percent with unequal diffusivities (Figure~\ref{fig:uneq}(c)); here the residual is the sum of all terms of \eqref{eq:energy-identity}, brought to one side and evaluated on the computed solution with $d\ee/dt$ replaced by its difference quotient in time (Appendix~\ref{app:scheme}), divided by the sum of their absolute values, and the smaller value with unequal diffusivities reflects the larger sum that normalizes it. With unequal diffusivities, $\mathcal A$ levels off at about five times the Joule heating $\Jou$ (Figure~\ref{fig:uneq}(d)). Since $\epsilon_+\ll\epsilon_-$ and the discharge is, as with equal diffusivities, a positive ion cloud, $\mathcal A$ nearly cancels against the dissipative term $\frac{\epsilon_++\epsilon_-}{2\epsilon_0}\norm q_2^2$: the two curves coincide to within the line width. Their difference $\frac1{\epsilon_0}\int_\Om(\epsilon_-n-\epsilon_+p)\,q\dx$, which is nonpositive under (H1), is positive here and, once the discharge is stationary, about two percent of $\Jou$. 
Hence, in this simulation, once the discharge is stationary the integrand of \eqref{eq:needed} satisfies
$\mathcal A\le\delta\|q\|_2^2+\beta\mathcal J$ with $\delta=0.95\,\frac{\epsilon_++\epsilon_-}{2\epsilon_0}$
and $\beta$ about a quarter: at $10$ ns, $\mathcal A$ and the dissipative term are both about $760$ W/m, their difference about $3$ W/m and $\Jou$ about $160$ W/m. Section~\ref{ssec:unequal} states the estimate that is missing.


\section{Discussion and open problems}\label{sec:remarks}

Theorem~\ref{thm:B} settles the question of \cite[Remark~1.2]{GHL26} under (H1) and (H2). We discuss what remains open without (H1) and where (H2) fails (Sections~\ref{ssec:unequal} and~\ref{ssec:touchdown}), which singularities Theorems~\ref{thm:A} and~\ref{thm:B} leave possible (Section~\ref{ssec:scenarios}), and what the theorems and the simulations suggest about the breakdown threshold (Section~\ref{ssec:threshold}).

\subsection{Unequal diffusivities}\label{ssec:unequal}
Under (H2), Proposition~\ref{prop:conditional} reduces global existence for $\epsilon_+\neq\epsilon_-$ to the a priori bound \eqref{eq:needed} on the ambipolar term $\mathcal A$ of \eqref{eq:ambipolar}, the only term of the balance law \eqref{eq:energy-identity} that the closure of Section~\ref{sec:global} does not absorb. We have not been able to prove \eqref{eq:needed}; we record what the simulation of Section~\ref{ssec:num-uneq} shows and why the obvious routes to it fail.

In that simulation the electrostatic energy and the number of carriers, the quantities of Theorem~\ref{thm:B}, differ by less than two percent from those with equal diffusivities. Once the discharge is stationary, the computed solution satisfies $\mathcal A\le\delta\norm q_2^2+\beta\Jou$ with $\beta<1$ for $\delta$ close to $\frac{\epsilon_++\epsilon_-}{2\epsilon_0}$, because $\mathcal A$ nearly cancels against the dissipative term $\frac{\epsilon_++\epsilon_-}{2\epsilon_0}\norm q_2^2$ (Figure~\ref{fig:uneq}(d)). What is missing is this estimate for every weak solution. In terms of the difference between $\mathcal A$ and the dissipative term, \eqref{eq:needed} requires, up to its last two terms, that
\[
\int_0^t\Big(\mathcal A-\frac{\epsilon_++\epsilon_-}{2\epsilon_0}\norm q_2^2\Big)\dr\le\beta\int_0^t\Jou\dr-\Big(\frac{\epsilon_++\epsilon_-}{2\epsilon_0}-\delta\Big)\int_0^t\norm q_2^2\dr
\]
with $\beta<1$ and $\delta<\frac{\epsilon_++\epsilon_-}{2\epsilon_0}$. Heuristically, the condition $\beta<1$ is borderline: for a quasi-neutral state, $p=n=s/2$, in which the charge current of Remark~\ref{rem:interpretation} vanishes and the applied field is neglected, the space-charge field is the ambipolar field $\nabla\psi=\frac{\epsilon_--\epsilon_+}{\mu_++\mu_-}\nabla\log s$; the dissipative term is then negligible, and \eqref{eq:weak-poisson} gives, up to terms with the data $p_D,n_D$, $\mathcal A=\frac{\epsilon_--\epsilon_+}{2}\int_\Om\nabla\psi\cdot\nabla s\dx=\frac{\mu_++\mu_-}{2}\int_\Om s\abs{\nabla\psi}^2\dx=\Jou$ to leading order.

Three routes to \eqref{eq:needed} fail. Bounding $\mathcal A$ by the $L^2$-norms of the densities, and these by the estimates of Section~\ref{sec:L2}, reintroduces a loop of the Riccati type \eqref{eq:GHL-riccati}, which the proof of Theorem~\ref{thm:B} avoids. By Young's inequality, a bound on $\int_0^t\norm s_2^2\dr$ on bounded time intervals would suffice for global existence, though not for bounds of the form \eqref{eq:growth}; in dimension two, an entropy with the weights $1/\mu_\pm$ controls this integral through the Fisher information it dissipates, but only up to a factor $1+m$ and with a source term that carries the logarithm of the densities, and we have not been able to close this argument. Perturbing Theorem~\ref{thm:B} does not give global existence for small $\abs{\epsilon_+-\epsilon_-}$ either, since the bounds \eqref{eq:growth} grow at least like $(1+t)e^{Ct}$; near a stable stationary solution of the problem with equal diffusivities, if there is one (Problem~\ref{conj:threshold}), the question becomes one of perturbation and linearized stability.

\subsection{Hypothesis \textup{(H2)} and sparking at touch-down}\label{ssec:touchdown}
The proof of Theorem~\ref{thm:B} uses (H2) twice: in the battery term of Corollary~\ref{cor:energy-ineq} and in \eqref{eq:source-Joule}. We keep (H2) as a hypothesis because (A2) provides only $\phi_D\in H^2(\Om)$, which does not imply $\nabla\phi_D\in L^\infty(\Om)$. Hypothesis (H2) holds when $\GN$ consists of planar pieces that do not meet one another and the electrodes are $C^{1,1}$ surfaces that carry constant potentials and meet $\GN$ at right angles: near each planar piece of $\GN$, the even reflection of $\phi_D$ across its plane extends it to a bounded harmonic function in a domain whose boundary there is $C^{1,1}$ and carries constant data, which belongs to $W^{2,p}$ near the boundary for every $p<\infty$ \cite[Theorem~2.4.2.5]{Gri85}, and $W^{2,p}(\Om)\hookrightarrow W^{1,\infty}(\Om)$ for $p>d$. In the geometry of \cite[Figure~1]{GHL26}, where the electrodes meet the artificial boundary at right angles, (H2) therefore holds as long as the electrodes are apart; Section~\ref{ssec:num-setting} uses this argument for the domain \eqref{eq:num-domain}.

Above the pull-in voltage (Section~\ref{ssec:background}), the deformable electrode approaches the fixed one, and the minimal gap width $w_{\min}(t)$ tends to zero at a finite time $T_*$ \cite{EGG10,Guo08,KLNT15}. Since $\phi_D$ increases by $V$ across the narrowest part of the gap, $\norm{\nabla\phi_D}_\infty\ge V/w_{\min}(t)\to\infty$ as $t\uparrow T_*$, so (H2) fails at touch-down. The electron current density injected at the cathode, of order $\mu_-n_DV/w_{\min}(t)$, becomes unbounded as well. Under (H1)--(H2), Theorem~\ref{thm:B} excludes a finite-time singularity of \eqref{eq:sys} in a fixed domain. With equal diffusivities, a singular description of sparking must therefore involve the degeneration of the geometry at touch-down, which lies outside \eqref{eq:sys}, the growth as $t\to\infty$ discussed in Section~\ref{ssec:threshold}, or the limit $\epsilon_\pm\to0$ of Section~\ref{ssec:scenarios}.

The bounds of Theorem~\ref{thm:B} do not depend on $\alpha_2$ (Remark~\ref{rem:constants}); the gap width at which the multiplication of electrons is largest does. If diffusion and space charge are neglected, the electrons crossing a gap of width $w$ in the uniform field $V/w$ are multiplied by $e^{w\alpha(V/w)}$, with the net Townsend coefficient $\alpha(E)=\alpha_1e^{-\alpha_2/E}-\eta_0$, and $w\,\alpha(V/w)=\alpha_1\,w\,e^{-\alpha_2w/V}-\eta_0\,w$ is largest, up to the attachment term, at $w_*=V/\alpha_2$, which is $10.8$ $\mu$m for air at atmospheric pressure and $300$ V, close to the $10$ $\mu$m of the dip in Section~\ref{sec:numerics}. As $w\to0$ the multiplication tends to one, since $\abs{w\,\alpha(V/w)}\le\gamma_*w$ by \eqref{eq:Gamma-bound}, whereas the vacuum field and the injected current density diverge.

Theorems~\ref{thm:A} and~\ref{thm:B} bear on a moving electrode insofar as the time scales separate: the electrode moves on the mechanical time scale, microseconds and longer \cite{PB02}, whereas the computed discharges at $150$ and $300$ V level off within ten nanoseconds (Figure~\ref{fig:voltages}(a) and~(b)). The coupled problem, in which the electrostatic force moves the electrode and the discharge modifies that force, is the natural next step.

\subsection{Which singularities remain possible}\label{ssec:scenarios}
Theorems~\ref{thm:A} and~\ref{thm:B} constrain three mechanisms by which \eqref{eq:sys} could form a finite-time singularity ($T_{\max}<\infty$): blow-up of the amplitudes $\norm{p}_\infty,\norm{n}_\infty$ at a point, self-similar collapse at a shrinking scale, and a gradient catastrophe at a steepening front. Under (H1)--(H2), Theorem~\ref{thm:B} excludes all three as causes of a finite $T_{\max}$; it does not bound $\norm p_\infty$, $\norm n_\infty$ or the gradients in $L^\infty$, so a loss of boundedness that keeps the norms in \eqref{eq:growth} finite is not addressed. Under (H2) with $\epsilon_+\neq\epsilon_-$, each requires \eqref{eq:needed} to fail (Proposition~\ref{prop:conditional}). By Theorem~\ref{thm:A}, for any $\epsilon_\pm$, each mechanism requires $\norm{p-n}_{\HDd}\notin L^{2\kappa}(0,T_{\max})$: quasi-neutral blow-up, with large densities and small charge, is impossible. We comment on the last two mechanisms.

\emph{Self-similar collapse.} Under the natural scaling of the drift--diffusion--Poisson system
the ionization source is of lower order.

\begin{proposition}[Scaling of the ionization source]\label{prop:scaling}
For $\lambda>0$ let $p_\lambda(x,t)=\lambda^2p(\lambda x,\lambda^2t)$, $n_\lambda(x,t)=\lambda^2n(\lambda x,\lambda^2t)$, $\phi_\lambda(x,t)=\phi(\lambda x,\lambda^2t)$ and $\bv_\lambda(x,t)=\lambda\bv(\lambda x,\lambda^2t)$. If $(\phi,p,n)$ solves \eqref{eq:sys} in $\Om\times(0,T)$ with the constants $(\alpha_1,\alpha_2,\eta_0)$ and the velocity $\bv$, then $(\phi_\lambda,p_\lambda,n_\lambda)$ solves \eqref{eq:sys} in $\lambda^{-1}\Om\times(0,\lambda^{-2}T)$ with the constants $(\lambda\alpha_1,\lambda\alpha_2,\lambda\eta_0)$ and the velocity $\bv_\lambda$.
\end{proposition}

\begin{proof}
The scaling multiplies the Poisson equation by $\lambda^2$ and the left-hand sides of the carrier equations by $\lambda^4$, since, with the unscaled quantities evaluated at $(\lambda x,\lambda^2t)$, $\partial_tp_\lambda=\lambda^4\partial_tp$, $\Delta p_\lambda=\lambda^4\Delta p$, $\nabla\cdot(p_\lambda\nabla\phi_\lambda)=\lambda^4\nabla\cdot(p\nabla\phi)$, $\nabla\cdot(p_\lambda\bv_\lambda)=\lambda^4\nabla\cdot(p\bv)$, $\Delta\phi_\lambda=\lambda^2\Delta\phi$ and $p_\lambda-n_\lambda=\lambda^2(p-n)$. The factor $n|\nabla\phi|$ of the source, by contrast, is multiplied only by $\lambda^3$: $n_\lambda|\nabla\phi_\lambda|=\lambda^3n|\nabla\phi|$. With the constants $(\lambda\alpha_1,\lambda\alpha_2,\lambda\eta_0)$ the net Townsend coefficient becomes $\lambda\alpha_1e^{-\lambda\alpha_2/|\nabla\phi_\lambda|}-\lambda\eta_0=\lambda\big(\alpha_1e^{-\alpha_2/|\nabla\phi|}-\eta_0\big)$, so that the source is multiplied by $\lambda^4$ as well.
\end{proof}

A self-similar collapse at a point would have densities of order $\ell^{-2}$ and field of order $\ell^{-1}$ on a spatial scale $\ell\propto(T_{\max}-t)^{1/2}$. Rescaling a neighbourhood of size $\ell$ to a fixed size $\ell_0$ is Proposition~\ref{prop:scaling} with $\lambda=\ell/\ell_0$ and multiplies the net Townsend coefficient by the factor $\ell/\ell_0$, which vanishes as $t\uparrow T_{\max}$: on scales below the ionization length $\alpha_1^{-1}$ the source is negligible against drift and diffusion, and what remains is the globally well-posed bipolar system of Remark~\ref{rem:freeenergy}. 
This suggests that ionization does not drive a self-similar collapse but acts at the fixed scale
$\alpha_1^{-1}$, where it produces Townsend growth. The argument is heuristic, and for unequal
diffusivities the ambipolar term $\mathcal A$ of \eqref{eq:ambipolar} remains the term that we cannot control.

\emph{Gradient catastrophe.} Under (H1)--(H2), the second bound in \eqref{eq:growth} keeps $\nabla p$ and $\nabla n$ square integrable on bounded time intervals, so that diffusion regularizes steep fronts at fixed $\epsilon$; but the proof of Theorem~\ref{thm:B} gives no bound uniform in $\epsilon$, and the theorem says nothing about the limit $\epsilon_\pm\to0$, in which the carrier equations become first-order transport equations coupled through the Poisson equation. Heuristically, planar ionization fronts, whose width is set by diffusion (Section~\ref{ssec:background}), sharpen in that limit. Whether, as $\epsilon_\pm\to0$, solutions of \eqref{eq:sys} develop layers of vanishing width with unbounded density gradients is open.

\subsection{The breakdown threshold}\label{ssec:threshold}
A natural mathematical counterpart of the breakdown voltage would be a threshold voltage below which solutions are global and above which they blow up in finite time. Theorem~\ref{thm:B} rules out such a threshold: under (H1)--(H2) every solution is global, at every voltage. A breakdown threshold, if there is one, must therefore separate two global behaviours, which we expect to be relaxation to a stationary state and unbounded growth as $t\to\infty$. The simulations of Section~\ref{ssec:num-voltage} show only the first: at $150$ and $300$ V the solution relaxes to what appears to be a stationary state, at $450$ V it overshoots and is still decreasing when the simulation stops, and at all three voltages the electrostatic energy and the $\HDd$-norm of the charge, which by Theorem~\ref{thm:A} must be unbounded near a finite-time singularity, level off or decrease (Figure~\ref{fig:voltages}(b) and~(c)). This leads to the following question.

\begin{problem}[Breakdown threshold]\label{conj:threshold}
Let Assumption~\ref{ass:main}, (H1) and (H2) hold, let the electrode potential be $V$ on the
deformable electrode and $0$ on the fixed one, and let the domain, the carrier data $p_D,n_D$ and a time-independent velocity field $\bv$ be fixed. Is there a threshold $V_*\in(0,\infty)$ such that for $0<V<V_*$ problem \eqref{eq:sys}--\eqref{eq:bc} has a stationary solution which is stable, in the sense that the solutions with initial data close to it in $L^2(\Om)$ converge to it in $L^2(\Om)$ as $t\to\infty$, whereas for $V>V_*$ there are initial data for which
$\|p(t)\|_2+\|n(t)\|_2\to\infty$ as $t\to\infty$?
\end{problem}

Two precedents motivate the question. For the corona model of \cite{Bud91,BW90} between coaxial cylinders, of the form \eqref{eq:sys} (Section~\ref{ssec:background}), the computed stationary solution was observed to lose stability at high voltage. In
the model of Degond, Lucquin-Desreux and Morrow, breakdown is the loss of stability of the
non-ionized gas at the sparking voltage, from which a branch of stationary discharges
bifurcates \cite{SS19,SS25,ST20}; this threshold exists without secondary emission
\cite{ST20}, in contrast to Townsend's criterion. In \eqref{eq:sys}--\eqref{eq:bc} the
electrode data (A3) exclude a non-ionized state, and $V_*$ concerns the stability of a
stationary discharge that carries space charge. The stationary problem is studied in \cite{HX26}; beyond the existence of stationary
solutions, the question involves a nonlinear stability problem whose
linearized operator is non-self-adjoint and nonlocal, and the passage from linear instability
to unbounded growth. For $\epsilon_+\ne\epsilon_-$, where global existence is open
(Section~\ref{ssec:unequal}), a threshold cannot yet be described as one between two global
behaviours.

\subsection*{Acknowledgements}
The author thanks Heiko Gimperlein and Andrew Lacey for the collaboration \cite{GHL26} which motivated this work.

\appendix

\section{The numerical scheme}\label{app:scheme}
We describe the scheme used in Section~\ref{sec:numerics}. The domain is the region \eqref{eq:num-domain}, with $\GD=\mathsf A\cup\mathsf B$ and $\GN=\mathsf C=\{|x|=r\}$; the gas is at rest, $\bv=0$, and the electrode data $p_D=n_D=n_b$ are constant. The scheme is a cell-centred finite volume method on the uniform Cartesian grid of $200\times150$ cells of size $h_x\times h_y=0.4\times0.2$ $\mu$m covering the rectangle $(-r,r)\times(0,w_{\max})$, restricted to the cells whose centres lie in $\Om$; the electrode $\mathsf A$ is represented by the staircase of faces separating these cells from the others. Dirichlet conditions enter through the faces on the electrodes, at the distance $h/2$ from the adjacent centres ($h=h_y$ on horizontal faces, $h=h_x$ on vertical faces), and the Neumann conditions as zero fluxes through the faces on $\mathsf C$. The harmonic extension $\phi_D$ is computed on the same grid with the same discrete Laplacian, so that $\psi=\phi-\phi_D$ vanishes on the electrodes exactly; the corners of the staircase enhance the vacuum field locally, and the largest cell value of $|\nabla\phi_D|$ is $1.35\,V/w_{\min}$, attained at the first step of the staircase next to the bottom of the dip.

\smallskip\noindent\emph{Fluxes.} For a carrier density $c$ with diffusivity $\epsilon_c$ and drift velocity $u_c\nabla\phi$ ($\epsilon_c=\epsilon_-$, $u_c=\mu_-$ for the electrons; $\epsilon_c=\epsilon_+$, $u_c=-\mu_+$ for the ions), the flux $-\epsilon_c\nabla c+c\,u_c\nabla\phi$ through the face $\sigma=K|K'$ between neighbouring cells whose centres are at the distance $h_\sigma$ is approximated by the Scharfetter--Gummel flux \cite{Kul95,SG69}
\begin{equation}\label{eq:SG}
\mathsf F_{K|K'}(c,\phi)=\frac{\epsilon_c}{h_\sigma}\Big(B(-z_\sigma)\,c_K-B(z_\sigma)\,c_{K'}\Big),\qquad
z_\sigma:=\frac{u_c\,(\phi_{K'}-\phi_K)}{\epsilon_c},\qquad B(z):=\frac{z}{e^{z}-1},
\end{equation}
which is the central difference of the diffusive flux for $z_\sigma=0$ and the upwind drift flux for $|z_\sigma|\to\infty$; on the faces on the electrodes the neighbour value is the Dirichlet datum and $h_\sigma=h/2$. Poisson operators $-\nabla\cdot(a\nabla\phi)$ with a cell-wise coefficient $a$ are discretized by the two-point flux $-a_\sigma(\phi_{K'}-\phi_K)/h_\sigma$ with the arithmetic mean $a_\sigma$ of the adjacent values. We write $|\sigma|$ for the length of the face $\sigma$ and $(\nabla\phi)_\sigma:=(\phi_{K'}-\phi_K)/h_\sigma$ for the face gradient.

\smallskip\noindent\emph{Time stepping.} Let $(\phi^k,p^k,n^k)$ approximate the solution at time $t_k$, $q^k=p^k-n^k$, $g^k=\mu_+p^k+\mu_-n^k$, and $\tau_k=t_{k+1}-t_k$. A step consists of three linear solves. For $\bv=0$ the charge satisfies $\partial_tq=\nabla\cdot(\epsilon_+\nabla p-\epsilon_-\nabla n)+\nabla\cdot(g\nabla\phi)$; inserting $q^{k+1}\approx q^k+\tau_k\partial_tq$ into $-\epsilon_0\Delta\phi^{k+1}=q^{k+1}$, with the last term implicit in $\phi$, gives the semi-implicit Poisson equation \cite{HK00,VHK94}
\begin{equation}\label{eq:semi-implicit}
-\nabla\cdot\big((\epsilon_0+\tau_k g^k)\nabla\phi^{k+1}\big)=q^k+\tau_k\,\nabla\cdot\big(\epsilon_+\nabla p^k-\epsilon_-\nabla n^k\big),
\end{equation}
with $\phi^{k+1}=\phi_D$ on the electrodes and $\bnu\cdot\nabla\phi^{k+1}=0$ on $\GN$; it removes the restriction of the time step to the dielectric relaxation time $\epsilon_0/\max g$ that an explicit coupling would impose. The carrier equations are then advanced by the backward Euler method with the fluxes \eqref{eq:SG} evaluated at $\phi^{k+1}$: the electron equation with the ionization rate $\Gamma(|\nabla\phi^{k+1}|_K)$ implicit, which the linearity of the source in $n$ allows, $|\nabla\phi|_K$ being the modulus of the central-difference gradient of the cell, and the ion equation with the source $\Gamma(|\nabla\phi^{k+1}|_K)\,n^{k+1}_K$ now known, made proportional to $p^{k+1}_K$ where $\Gamma<0$ and it is a sink. The source enters the two equations with the same amount wherever $\Gamma\ge0$, so that it cancels in the discrete charge balance as in \eqref{eq:charge-eq}, and the matrices of the two solves are M-matrices for $\tau_k\Gamma<1$, so that the densities stay nonnegative; in all simulations the minimum of $n$ was $0.956\,n_b$. The time step is $\tau_k=\min\{1\ \mathrm{ps},\,0.7/\max_K\Gamma(|\nabla\phi^{k}|_K)\}$, which gives $\tau_k$ between $0.25$ and $1$ ps; the simulations cover $10$ ns, or $5$ ns at $450$ V, in $10^4$ to $2\cdot10^4$ steps.

\smallskip\noindent\emph{Evaluation and accuracy.} The quantities of Section~\ref{sec:numerics} are evaluated from the discrete fields, with $\psi=\phi-\phi_D$: $\ee=\frac{\epsilon_0}{2}\sum_\sigma|\sigma|h_\sigma(\nabla\psi)_\sigma^2$ over all faces including those on the electrodes, $\norm{q}_{\HDd}=\epsilon_0\norm{\nabla\psi}_2$, and the $L^1$- and $L^2$-norms as sums over the cells. The Joule heating $\Jou$ and the battery term are evaluated in the same way with the face gradients and the arithmetic means of $g$ over the cells adjacent to each face, and $d\ee/dt$ is replaced by the difference quotient in time. The relative residual of the balance law \eqref{eq:energy-identity} defined in Section~\ref{ssec:num-uneq} and shown in Figure~\ref{fig:uneq}(c) comes from the numerical diffusion of the upwinding in \eqref{eq:SG}: at $300$ V with equal diffusivities it is $1.5$ percent, and $2.2$ percent on a grid coarser by the factor $1.5$ (medians over the second half of the simulation), consistent with the first-order accuracy of the upwinding; with the diffusivities of air it is $0.4$ percent, the sum of terms that normalizes it being then dominated by $\mathcal A$ and the dissipative term. At $300$ V, the coarser grid reproduces the electrostatic energy, the number of carriers, the $\HDd$-norm of the charge and the maximal ion density at $10$ ns to within $10$ to $21$ percent, and the electrode fields to within $9$ percent; the leading error is the overestimate of the avalanche amplification by the upwind treatment of the source, by the factor $(1-\alpha h_y)^{-1}$ in place of $e^{\alpha h_y}$ per cell along the drift of the electrons, $\alpha$ being the net Townsend coefficient; at the mean field $V/w_{\min}$ of $300$ V this amounts to about $24$ percent across the dip on the grid of Section~\ref{sec:numerics} and to $39$ percent on the coarser grid, so that the error is of first order in the mesh size.

\end{document}